\input amsppt.sty
\loadbold
\TagsOnRight
\hsize 30pc
\vsize 47pc
\magnification 1100
\redefine\g{\goth g}
\redefine\su{\goth{su}}

\redefine\SU{ {\text{SU}}}
\redefine\SO{ {\text{SO}}}
\redefine\Sp{{\text{Sp}}}
\redefine\a{\goth a}
\redefine\B{{\Cal B}}
\redefine\D{{\Cal D}}
\define\F{{\Cal F}}
\define\Hom{\operatorname{Hom}}

\redefine\l{\goth l}
\redefine\k{\goth k}
\redefine\h{\goth h}
\redefine\m{\goth m}
\redefine\n{\goth n}
\redefine\p{\goth p}
\redefine\t{\goth t}
\redefine\z{\goth z}

\define\r{{\goth r}}
\define\q{{\goth q}}
\redefine\C{\Bbb C}
\redefine\R{\Bbb R}
\define\Z{\Bbb Z}

\redefine\span{\text{span}}

\define\={\overset\text{def}\to=}

\define\RP{\R P}
\define\CP{\C P}

\topmatter
\title The Ricci tensor of an almost homogenous K\"ahler manifold \endtitle
\author  Andrea Spiro \endauthor 
\address
\newline
Andrea Spiro\newline
Mathematics and Physics Department\newline
 Via Madonna delle Carceri\newline
I- 62032 Camerino (Macerata)\newline
ITALY\newline
\newline
 \endaddress 
\email 
spiro\@campus.unicam.it
\endemail 
\keywords K\"ahler-Einstein metrics, cohomogeneity one 
actions\endkeywords
\subjclass  53C55, 53C25, 57S15\endsubjclass
\thanks \endthanks 

 \abstract We determine an explicit expression for the Ricci tensor of 
a K-manifold, that 
is of a compact K\"ahler 
manifold $M$ with vanishing first Betti number, on which a semisimple group 
$G$ of 
biholomorphic isometries acts with an orbit of codimension one. We also prove that 
the K\"ahler form $\omega$ and the Ricci form $\rho$ of $M$  are uniquely determined
by two special curves with values in 
$\g = Lie(G)$, say $Z_\omega, Z_\rho: \R \to \g = Lie(\g)$ and we show how 
 the curve $Z_\rho$ is determined by the curve $Z_\omega$.\par
These results are used in another work with F. Podest\`a, where  
new examples of non-homogeneous
compact K\"ahler-Einstein manifolds with positive first Chern class are 
constructed.
\endabstract

\leftheadtext{\smc Ricci Tensor of an almost K\"ahler homogeneous manifold}
\rightheadtext{\smc A. Spiro}

\endtopmatter
\document

\subhead 1. Introduction
\endsubhead
\bigskip
The objects of our study are the 
so-called {\it K-manifolds\/}, that  is  K\"ahler manifolds $(M,J,g)$ with $b_1(M) = 0$
and which are acted on by a group $G$
of biholomorphic isometries, 
with regular orbits of codimension one. 
Note that  since $M$ is compact and $G$ has orbits of codimension one, 
the complexified group $G^\C$ acts naturally on $M$
as a group of biholomorphic transformations, with an
open and dense orbit.  According to a terminology introduced by 
A. Huckleberry and D. Snow in \cite{HS},  $M$ is  {\it almost-homogeneous\/}
with respect to the $G^\C$-action. By the results in \cite{HS}, 
the subset $S\subset M$ of singular points for the $G^\C$-action is either
connected or 
with exactly two connected components. If the first case occurs,
we will say that $M$ is  a {\it non-standard K-manifold\/}; 
we will call it {\it standard K-manifold\/} in the 
other case. 
\par
The aim of this paper is furnish  an explicit expression for  
the Ricci curvature tensor of a K-manifold, to
be used 
for constructing (and possibly classify) new families 
of examples of non-homogeneous 
 K-manifold with special curvature conditions. A successful application 
of our results is given 
in \cite{PS1}, where several new examples of non-homogeneous 
compact K\"ahler-Einstein 
manifolds with positive first Chern class are found.\par
\medskip
Note that explicit expressions for the Ricci tensor of 
standard K-manifolds can be  found also 
in \cite{Sa}, \cite{KS}, \cite{PS} and \cite{DW}. However our results
can be applied to any kind of  K-manifold and hence they turn out to be 
 particularly useful 
for the non-standard cases (at this regard, see also \cite{CG}). 
They can be resumed in 
the following three facts.\par
\medskip
Let $\g$ be the Lie algebra of the compact group $G$ acting on
the K-manifold $(M, J,g)$ with at least one orbit of codimension one. 
By a result of \cite{PS1}, 
we may always assume that $G$ is semisimple. Let also $\B$ be
the Cartan-Killing form
of $\g$. Then  for any   $x$ in the regular point set $M_{\text{reg}}$, 
one can consider the following $\B$-orthogonal decomposition of $\g$:
$$\g = \l + \R Z + \m\ , \tag 1.1$$
where $\l = \g_x$ is the isotropy subalgebra, $\R Z + \m$ is naturally 
identified with the
tangent space $T_o(G/L) \simeq T_x(G\cdot x)$ of the $G$-orbit $G/L = G\cdot x$, and $\m$ 
is naturally 
identified  with the holomorphic
 subspace $\m \simeq \D_x$ 
$$\D_x = \{\ v\in T_x(G\cdot x) \ :\ Jv \in  T_x(G\cdot x)\ \}\ .\tag1.2$$
 Notice that for any  point $x \in M_{\text{reg}}$ 
the $\B$-orthogonal decomposition (1.1)  is uniquely given; 
on the other hand, two distinct points
$x, x' \in M_{\text{reg}}$ may determine
 two distinct decompositions of type
(1.1).\par
\smallskip
 Our first  result consists 
in proving that  any  K-manifold
has a family $\Cal O$ of smooth
 curves $\eta: \R \to M$  of the form
$$\eta_t = \exp(i t Z)\cdot x_o\ ,$$
where $Z \in \g$,  $x_o\in M$ is a regular point for the $G^\C$-action and 
the following properties are satisfied:
\roster
\item $\eta_t$ intersects any regular $G$-orbit;
\item for any  point $\eta_t \in M_{\text{reg}}$,  the tangent vector
$\eta'_t$
is transversal to the regular orbit $G\cdot \eta_t$; 
\item
any element $g\in G$ which belongs to a 
stabilizer $G_{\eta_t}$, with $\eta_t\in M_{\text{reg}}$,
fixes pointwise the whole curve $\eta$; in particular, all regular
orbits $G\cdot \eta_t$  are equivalent to the same
homogeneous space $G/L$;
\item the decompositions (1.1) associated with the points   $\eta_t \in M_{\text{reg}}$
do not depend on $t$;
\item there exists a basis $\{f_1, \dots f_n\}$ for $\m$ such that for any 
$\eta_t \in M_{\text{reg}}$
the complex structure $J_t: \m \to \m$,  induced  by the complex structure 
of $T_{\eta_t}M$,  is of the following form:
$$J_t f_{2j} =  \lambda_j(t) f_{2j+1}\ ,\qquad J_t f_{2j+1} =  - \frac{1}{\lambda_j(t)}f_{2j}\ ;
\tag1.3$$
where the function $\lambda_j(t)$ is either one of the functions
  $- \tanh(t)$,  $- \tanh(2t)$, $-\coth(t)$ and 
 $- \coth(2t)$ or it is identically equal to  $1$.
\endroster
We call any such  curve  an {\it optimal transversal curve\/};  the 
 basis for  $\R Z + \m \subset \g$ given by 
$(Z, f_1, \dots, f_{2n-1})$, 
where the  $f_i$'s  verify (1.3),
 is called {\it optimal basis associated with $\eta$\/}. An
explicit description of the optimal basis 
for any given semisimple Lie group $G$ is given 
in 
\S 3.\par
 Notice that the family $\Cal O$ of
optimal transversal curves  depends only on the 
action of the Lie group $G$. In particular 
it is totally independent on the choice of the $G$-invariant K\"ahler
metric $g$. At the same time, the Killing fields, associated with the elements 
of an optimal basis, determine a 1-parameter family 
of holomorphic  frames at the points $\eta_t \in M_{\text{reg}}$, which are
orthogonal  w.r.t. at least one  $G$-invariant K\"ahler metric $g$. 
It is also proved that, for all K-manifold $M$ which do not 
belong to a special class of non-standard
K-manifold, those holomorphic frames are 
orthogonal w.r.t. {\it any} $G$-invariant K\"ahler metric $g$ on $M$
 (see Corollary 4.2 for details).
From these remarks and the fact that  $\eta'_t = J\hat Z_{\eta_t}$, 
where $Z$ is the first element of any  optimal basis, it may be inferred that any 
curve $\eta \in \Cal O$ is  a reparameterization of a 
normal geodesics of some (in most cases, {\it any\/})
$G$-invariant K\"ahler metric  on $M$.\par
\medskip
Our second main  result
is the following. Let $\eta$ be an optimal transversal 
curve of a  K-manifold, $\g = \l + \R Z + \m$ the decomposition (1.1) 
associated with the regular points  $\eta_t\in M_{\text{reg}}$ and  let
$\omega$ and $\rho$ be the K\"ahler form and the Ricci form, 
respectively, associated with a 
given $G$-invariant K\"ahler metric $g$ on $(M,J)$. \par
By a slight modification of  arguments used in \cite{PS}, we show 
that 
there exist two smooth curves 
$$Z_\omega, Z_\rho : \R \to C_\g(\l) = \z(\l) + \a\ ,\qquad  
\a = C_\g(\l) \cap (\R Z + \m)\ ,\tag 1.4$$ 
satisfying 
the following properties (here 
$\z(\l)$ denotes the center of $\l$ and 
$C_\g(\l)$ denotes the centralizer of $\l$ in $\g$): 
for any  $\eta_t\in M_{\text{reg}}$
and any two element $X, Y \in \g$, with associated 
Killing fields $\hat X$ and 
$\hat Y$, 
$$\omega_{\eta_t}(\hat X, 
\hat Y) = \B(Z_\omega(t), [X, Y])\ ,\quad
\rho_{\eta_t}(\hat X, 
\hat Y) = \B(Z_\rho(t), [X, Y])\ .\tag 1.5$$
We call such curves $Z_\omega(t)$ and $Z_\rho(t)$  {\it the algebraic representatives of 
$\omega$ and  $\rho$ along $\eta$\/}.  It is clear that the algebraic 
representatives determine uniquely the restrictions of $\omega$ and $\rho$
to the tangent spaces of the regular orbits. 
But the following Proposition establishes  a result which is somehow stronger.\par
Before stating the proposition, we recall that 
in \cite{PS} the following fact was established: if
$\g = \l + \R Z + \m$ is a decomposition of the form (1.1), then 
the  subalgebra  
$\a = C_\g(\l) \cap (\R Z + \m)$
is   either  1-dimensional or 3-dimensional 
and isomorphic with $\goth{su}_2$. By virtue of this dichotomy, the two cases
considered in the following proposition are  all possible cases.\par
\proclaim{Proposition 1.1} 
Let  $\eta_t$ be an optimal transversal curve
of a K-manifold $(M, J, g)$ acted on by the compact semisimple Lie group $G$
and  let $\g = \l + \R Z + \m$ be the decomposition of the form (1.1) determined 
by the points $\eta_t \in M_{\text{reg}}$. Let also 
$Z: \R \to C_\g(\l) = \z(\l) +  \a$ be
the algebraic representative of the K\"ahler form $\omega$ or 
of the Ricci form $\rho$. Then:\par
\roster
\item if  $\a$ is 1-dimensional, then it is of the form $\a = \R Z$ and
 there exists an element $I \in \z(\l)$ and 
a smooth function ${f}: \R \to \R$ so that
$$Z(t) = f(t) Z +  I\ ;\tag 1.6$$
\item if $\a$ is 3-dimensional, then it is  of the form $\a =  \goth{su_2} = \R Z+ \R X + \R Y$, 
with $[Z, X] = Y$ and $[X,Y] = Z$ and 
 there exists an element $I \in \z(\l)$, a real number $C$ and 
a smooth function $f: \R \to \R$ so that
$$Z(t) = f(t) Z_\D + \frac{C}{\cosh(t)}
X +  I\ .\tag 1.7$$
\endroster
Conversely,  if $Z: \R \to C_\g(\l)$ 
is a curve in $C_\g(\l)$ of the form (1.6) or 
(1.7), then there exists a 
unique closed $J$-invariant, $G$-invariant 2-form 
$\varpi$  on the set of regular points
$M_{\text{reg}}$, having $Z(t)$ as  algebraic representative.\par
 In particular,
the K\"ahler form $\omega$ and the Ricci form $\rho$ are uniquely determined by 
their  algebraic representatives.
\endproclaim
\bigskip
Using (1.5), Proposition 1.1 and some basic properties
 of the decomposition 
$\g = \l + \R Z + \m$ (see \S 5), it can be shown that
the algebraic representatives $Z_\omega(t)$ and $Z_\rho(t)$ 
are uniquely determined by the values 
 $\omega_{\eta_t}(\hat X, 
J \hat X) = \B(Z_\omega(t), [X, J_t X])$ and 
$\rho_{\eta_t}(\hat X, 
J \hat X) = \B(Z_\rho(t), [X, J_t X])$, where $X\in \m$ and $J_t$ is the complex structure 
on $\m$ induced by the complex structure of the tangent space $T_{\eta_t} M$. \par
\medskip
Here comes our third main result. It 
consists in Theorem 5.1 and Proposition 5.2, where we give the explicit  expression for 
the value $r_{\eta_t}(X,X) = \rho_{\eta_t}(\hat X, 
J \hat X)$  for any $X\in \m$, only in terms
of the algebraic representative $Z_\omega(t)$  and 
of the Lie brackets between  $X$ and the 
 elements of the optimal basis in  $\g$. By the previous discussion, 
this result furnishes a way to write down explicitly the Ricci tensor 
of the K\"ahler metric associated 
with $Z_\omega(t)$. \par
\bigskip
\remark{Acknowledgement} Many crucial ideas for this paper
are  the natural fruit of the uncountable discussions that 
Fabio Podest\`a and the author had  since they began 
working on cohomogeneity one K\"ahler-Einstein manifolds. 
It is fair to say that most of the credits should be shared with Fabio.
\endremark
\medskip
\remark{ Notation} Throughout the paper, if $G$ is a Lie group acting isometrically on a Riemannian 
manifold $M$ and if $X\in \g = Lie(G)$, 
we will  adopt the symbol $\hat X$  to denote the 
 Killing vector field on $M$ corresponding to $X$.\par
The  Lie algebra of a Lie
group  
will be always denoted  by the corresponding gothic letter. For a group $G$
and a Lie algebra $\g$, 
$Z(G)$ and $\z(\g)$  denote the center of $G$
and of $\g$, respectively. For any subset $A$ of a group $G$
or of a Lie algebra $\g$, $C_G(A)$ and $C_\g(A)$ are
the centralizer of $A$ in $G$ and $\g$, respectively. \par
Finally, for any subspace $\n \subset \g$ of a semisimple Lie algebra $\g$, 
the symbol $\n^\perp$ denotes the 
orthogonal complement of $\n$ in $\g$ w.r.t. the Cartan-Killing form $\B$.
\endremark
\bigskip
\bigskip
\subhead 2.  Fundamentals of K-manifolds 
\endsubhead
\bigskip
\subsubhead 2.1 K-manifolds, KO-manifolds and KE-manifolds
\endsubsubhead
\medskip
A {\it K-manifold\/} is a pair formed by a compact K\"ahler manifold $(M,J,g)$ 
and  a compact semisimple Lie group $G$
acting almost effectively and isometrically
(hence  biholomorphically)
on $M$, such that:
\roster
\item"i)"  $b_1(M)= 0$;
\item"ii)" $G$ acts of {\it cohomogeneity 
one with respect to the action of $G$\/}, i.e. the regular $G$-orbits 
are 
of codimension one in $M$.
\endroster
In this paper,  $(M, J, g)$ will always denote a K-manifold of dimension $2n$, 
acted on by 
the compact semisimple Lie group $G$. We will denote by
$\omega(\cdot, \cdot)  = g(\cdot, J \cdot)$ the K\"ahler fundamental form and by 
 $\rho = r(\cdot, J\cdot)$
the Ricci form of $M$.\par
\smallskip
For the general properties of  cohomogeneity one manifolds and of  K-manifolds,
see e.g. \cite{AA}, \cite{AA1}, \cite{Br}, \cite{HS}, \cite{PS}. 
Here we  only recall  some  properties, 
which will be used in the paper.
\smallskip
If $p\in M$ is a  regular point, let us denote by $L = G_p$
the corresponding isotropy subgroup. 
Since $M$ is orientable, every 
regular orbit $G\cdot p$ is orientable. Hence
we may consider a unit normal vector field $\xi$, defined 
on the subset of regular points $M_{\text{reg}}$, which is orthogonal 
to any regular orbit. It is known (see \cite{AA1})
that any integral curve of  $\xi$ 
is a geodesic. 
Any such geodesic is  usually called {\it normal geodesic\/}. \par
A normal geodesic $\gamma$ through a point $p$
verifies the following properties: 
it intersects any $G$-orbit orthogonally; 
the isotropy subalgebra $G_{\gamma_t}$ at a regular point $\gamma_t$ is 
always equal
$G_p = L$ (see e.g. \cite{AA}, \cite{AA1}). We formalize  
these two facts 
in the following definition.\par
We  call 
{\it nice transversal curve through a  point $p\in M_{\text{reg}}$\/} any  
curve $\eta : \R \to M$ with $p\in \eta(\R)$ and such that:
\roster
\item"i)" it intersects any regular orbit;
\item"ii)" for any $\eta_t \in M_{\text{reg}}$
$$\eta'_t \notin T_{\eta_t}(G\cdot {\eta_t})\ ;\tag 2.1$$
\item"iii)" for any $\eta_t \in M_{\text{reg}}$,  $G_{\eta_t} = L = G_p$.
\endroster
\bigskip
The following  property of K-manifold has been proved in \cite{PS}.\par
\proclaim{Proposition 2.1}  Let $(M, J, g)$ be a K-manifold acted on by the 
compact semisimple Lie group $G$. Let also $p\in M_{\text{reg}}$ and 
$L = G_p$ the isotropy subgroup at $p$.
Then:
\roster
\item there exists an element $Z$ (determined up to scaling) so that 
$$\R Z \in C_\g(\l)\cap \l^\perp\ , \quad C_\g(\l + \R Z) = \z(\l) + \R Z\ ;\tag 2.2$$
in particular, the connected subgroup $K \subset G$ with subalgebra $\k = \l + \R Z$
is the isotropy subgroup of a flag manifold $F = G/K$;
\item the dimension of  $\a = C_\g(\l) \cap \l^\perp $ is either $1$ or $3$; 
in case $\dim_\R \a = 3$, 
then $\a$ is a subalgebra isomorphic to $\goth{su}_2$
and there exists a Cartan subalgebra 
  $\t^\C \subset \l^\C + \a^\C \subset \g^\C$ so that 
$\a^\C = \C H_\alpha + \C E_\alpha + \C E_{-\alpha}$ for some 
root $\alpha$ of the root system of $(\g^\C, \t^\C)$.
\endroster
\endproclaim
\bigskip
Note that if  for some regular point $p$ we have that 
$\dim_\R \a = 1$  (resp. $\dim_\R \a = 3$), then the same
occurs 
at any other regular point. Therefore we may consider 
the following  definition.\par
\bigskip
\definition{Definition 2.2} Let $(M, J,g)$ be a K-manifold and 
$L = G_p$ the isotropy subgroup of a regular point $p$. We say that 
$M$ is a K-manifold {\it with ordinary action\/} (or shortly, {\it KO-manifold\/})
if $\dim_\R \a = \dim_\R(C_\g(\l) \cap \l^\perp) = 1$. \par
In all other cases, we say that $M$ is {\it with extra-ordinary action\/} (or, shortly, 
{\it KE-manifold\/}). 
\enddefinition
\bigskip
Another  useful property of K-manifolds is the following. 
It can be proved that any K-manifold admits exactly two singular orbits,
at least one of which is complex
(see \cite{PS1}).  By the results in \cite{HS}, it also 
follows that if $M$ is a K-manifold
whose singular orbits  are both complex, then $M$ admits a $G$-equivariant
blow-up $\tilde M$ along the complex singular orbits, which is still a K-manifold and 
admits a holomorphic
fibration over a flag manifold $G/K = G^\C/P$, with standard fiber
equal to $\C P^1$. \par
Several other important facts are related to the 
existence (or non-existence) of
two singular complex orbits (see \cite{PS1}
for a review of these properties). For this reason, it is convenient to introduce
the following definition.\par
\bigskip
\definition{Definition 2.3} We say that a K-manifold $M$, acted on by a
compact semisimple group $G$ with cohomogeneity one, is {\it standard\/} if the 
action of  $G$  has two 
singular complex orbits. We call it {\it non-standard\/} in all other cases.
\enddefinition
\bigskip
\subsubhead 2.2  The  CR structure 
of the regular orbits of a K-manifold
\endsubsubhead
\par
\bigskip
A {\it CR structure of codimension $r$\/} on  
a manifold $N$ is a pair $(\D, J)$ formed by a distribution $\D\subset TN$ of 
codimension $r$   and a smooth family 
$J$ of complex structures $J_x:
\D_x \to \D_x$ on the spaces  of the distribution. \par
A CR structure $(\D, J)$ is called {\it integrable\/}
if the distribution $\D^{10} \subset T^\C N$,  given by the
$J$-eigenspaces $\D^{10}_x \subset \D^\C_x$ corresponding to the 
eigenvalue $+i$, verifies
$$[\D^{10}, \D^{10}] \subset \D^{10}\ .$$
Note that  a complex structure $J$ on manifold $N$ may be 
always considered
as an integrable CR structure of codimension zero.\par 
\smallskip
A  smooth map $\phi: N \to N'$ between two CR manifolds $(N, \D, J)$
and $(N', \D', J')$ is  called
{\it CR  map\/} (or {\it holomorphic map\/}) if:
\roster
\item"a)" $\phi_*(\D) \subset \D'$; 
\item"b)" for any $x\in N$, 
$\phi_* \circ J_x = J'_{\phi(x)} \circ \phi_*|_{\D_x}$.  
\endroster
A {\it CR transformation\/} of $(N, \D, J)$ is 
a  diffeomorphism $\phi: N \to N$ which is also a CR map.\par
\medskip
Any codimension one submanifold $N \subset M$ of a complex manifold
$(M, J)$  is naturally endowed with an integrable  CR structure
 of codimension one $(\D, J)$, which is called {\it induced CR structure\/}; 
 it is defined by
$$\D_x = \{\ v\in T_xN\ :
\ J v\in T_xN\ \}\ \qquad \ \ J_x = J|_{\D_x}\ .$$\par
\medskip
It is clear that any regular orbit $G/L = G\cdot x \in M$
of a K-manifold $(M, J, g)$ has an induced CR structure $(\D,J)$, which 
is invariant under the transitive action of $G$. For this reason, several 
facts on the global structure of the regular orbits of a K-manifolds 
can be detected using what is known on compact homogeneous CR manifolds
(see e.g. \cite{AHR} and \cite{AS}).  \par
Here, we recall some of those facts, 
which will turn out to be crucial in the next sections.\par
\medskip
Let $(G/L, \D, J)$ be a  homogeneous CR manifold of a compact semisimple
Lie group $G$, with an integrable  CR structure $(\D, J)$ of codimension one.
If we consider the  $\B$-orthogonal decomposition 
$\g = \l + \n$, where $\l = Lie(L)$, then the orthogonal complement
$\n$ is naturally identifiable with 
the  tangent space $T_{o}(G/L)$, $o = eL$,  
by means  of the map 
$$\phi : \n \to T_{o}(G/L)\ ,\qquad \phi(X) = \hat X|_{o}\ .$$
If we denote by $\m$ the subspace
$$\m = \phi^{-1}(\D_{o}) \subset \n\ ,$$
 we get the 
following  orthogonal decomposition of $\g$:
$$\g = \l + \n = \l + \R Z_{\D} + \m\ .\tag 2.3$$
where $Z_\D  \in (\l + \m)^\perp$. Since the decomposition 
is $\operatorname{ad}_\l$-invariant, it follows that  $Z_\D
\in C_\g(\l)$. \par
\smallskip
Using again the identification map $\phi: \n \to T_{o}(G/L)$, we may consider the complex
structure
$$J: \m \to \m\ ,\qquad J  \= \phi^*(J_{o})\ .\tag 2.4$$
Note that   $J$ is uniquely determined by the 
 direct sum decomposition 
$$\m^\C = \m^{10} + \m^{01}\ ,\quad \m^{01} 
= \overline{\m^{10}}\ ,\tag 2.5$$
 where $\m^{10}$ and $\m^{01}$ are the 
$J$-eigenspaces  with  eigenvalues $+i$ and $-i$, respectively.\par
\smallskip
In all the following, (2.3) will be called 
{\it the structural decomposition of $\g$ associated with $\D$\/}; the subspace 
$\m^{10} \subset \m^\C$ (respectively, $\m^{01}
= \overline{\m^{10}}$) given (2.5) will be called {\it the holomorphic\/} (resp. 
{\it anti-holomorphic\/})
{\it subspace associated with $(\D,J)$\/}.\par
\bigskip
We recall that   a $G$-invariant CR structure $(\D, J)$ on 
$G/L$ is integrable if and only if the associated holomorphic subspace
$\m^{10} \subset \m^\C$ is so that 
$$\l^\C + \m^{10} \quad \text{is a subalgebra of }\quad \g^\C\ .\tag 2.6$$
\par
\bigskip
We  now need to introduce a few concepts which are quite helpful in describing 
the structure of a generic  compact homogeneous CR manifold.\par
\bigskip
\definition{Definition 2.4} Let $N = G/L$ be a homogeneous manifold of a 
compact semisimple Lie group $G$ and  $(\D,J)$  
a $G$-invariant, integrable CR structure  of codimension one on $N$. \par
We say that a  CR manifold $(N = G/L, \D, J)$ is a {\it 
Morimoto-Nagano space\/} 
if either  $G/L= S^{2n-1}$, $n>1$, endowed with the 
standard CR structure of $S^{2n-1} \subset \C P^n$,
or there exists a subgroup $H \subset G$ so that: 
\roster
\item"a)"
$G/H$ is   a compact rank
one symmetric space (i.e. $\R P^n = \SO_{n+1}/\SO_{n}\cdot \Z_2$,
$S^n = \SO_{n+1}/\SO_{n}$,  $\C P^n
= \SU_{n+1}/\SU_n$, $\Bbb H P^n = \Sp_{n+1}/\Sp_n$ or $\Bbb O P^2 = \text{F}_4/Spin_9$); 
\item"b)"  $G/L$ is 
a sphere bundle $S(G/H) \subset T (G/H)$ in the tangent space of $G/H$;
\item"c)" $(\D, J)$ is the 
CR structure induced  on $G/L = S(G/H)$ by the $G$-invariant
complex structure of $T(G/H) \cong G^\C/H^\C$.
\endroster
If a  Morimoto-Nagano space is    $G$-equivalent to a sphere $S^{2n-1}$
we call it {\it trivial\/}; 
we  call it  {\it  non-trivial\/} in all other cases.\par
A $G$-equivariant holomorphic fibering
$$\pi: N = G/L \to \F = G/Q$$
of $(N, \D, J)$ onto a non-trivial 
flag manifold $(\F = G/Q, J_\F)$  with invariant complex structure $J_\F$,
 is called {\it CRF fibration\/}.
A CRF fibration $\pi: G/L \to G/Q$ is called {\it nice\/} if the standard fiber
is a non-trivial Morimoto-Nagano space;  it 
is called {\it very nice\/} if it is nice and there exists no other 
nice CRF fibration $\pi': G/L \to G/Q$ with  standard fibers of smaller dimension.
\enddefinition
\bigskip
\medskip
The following Proposition  gives 
necessary and sufficient conditions for the existence of a CRF fibration. The proof
 can be found in  \cite{AS}.\par
\bigskip
\proclaim{Proposition 2.5} Let $G/L$ be  homogeneous CR manifold 
of a compact semisimple Lie group $G$, with an integrable, 
codimension one $G$-invariant CR structure $(\D, J)$. Let also
$\g = \l + \R Z_\D + \m$ be
 the   structural decomposition of $\g$ and  $\m^{10}$
the holomorphic subspace, associated with $(\D, J)$.\par
Then $G/L$  admits a non-trivial
CRF fibration if and only if there exists a proper parabolic 
subalgebra $\p = \r + \n \subsetneq \g^\C$ (here $\r$
is a reductive part and $\n$ the nilradical of $\p$) such that:
$$a)\ \r = (\p \cap \g)^\C\ ;\qquad b)\ \l^\C + \m^{01}\subset \p\ ;\qquad
c)\ \l^\C \subsetneq \r\ .$$
 In this 
case, $G/L$ admits a CRF fibration with basis $G/Q = G^\C/P$, where $Q$ is the connected 
subgroup generated by $\q = \r \cap \g$ and $P$ is the parabolic subgroup 
of $G^\C$ with Lie algebra $\p$. 
\endproclaim
\bigskip
Let us go back to the  regular orbits  of a K-manifold $(M, J,g)$
acted on by the compact semisimple group $G$. We already 
pointed out that  each regular orbit $(G/L = G\cdot x, \D, J)$, endowed with the 
induced CR structure $(\D, J)$, is 
a compact homogeneous CR manifold. In the statement of the following Theorem 
we collect the main results on the 
one-parameter family of compact homogeneous CR manifolds given by the regular orbits of 
a K-manifold, which is a direct consequence Th. 3.1 in \cite{PS1} (see also \cite{HS} and 
\cite{PS} Th.2.4). \par
\bigskip
\proclaim{Theorem  2.6}Ê Let $(M, J, g)$ be a K-manifold acted on by 
the compact semisimple Lie group $G$. \par
\roster
\item If $M$ is  standard,  then there exists 
a flag manifold  $(G/K, J_o)$ with a $G$-invariant complex structure $J_o$,
such that  any 
regular orbit $(G\cdot x = G/L$, $\D, J)$ of $M$ admits  
a CRF-fibration $\pi: (G/L, \D, J) \to (G/K, J_o)$ 
onto $(G/K, J_o)$
with standard fiber $S^1$.
\item If $M$ is  non-standard, then there exists 
a flag manifold  $(G/K, J_o)$ with a $G$-invariant complex structure $J_o$
such that any regular orbit $(G/L = G\cdot x, \D, J)$ admits a very 
nice CRF fibration 
$\pi: (G/L, \D, J) \to (G/K, J_o)$ 
where the standard fiber $K/L$ is a non-trivial Morimoto-Nagano space of 
dimension $\dim K/L \geq 3$.
\endroster
Furthermore, if the last case occurs, then
the fiber $K/L$ of the  CRF fibration $\pi: (G/L, \D, J) \to (G/K, J_o)$ 
has dimension  3 if and only if $M$ is a non-standard
KE-manifold and 
$K/L$ is  either $S(\RP^2) \subset T(\RP^2) = \CP^2\setminus\{\ [z]\ : {}^tz\cdot z = 0\}$
or $S(\CP^1) \subset T(\CP^1)  = \CP^1\times \CP^1\setminus \{\ [z] = [w]\ \}$. 
\endproclaim
\bigskip
\bigskip
\subhead 3. The optimal transversal curves of a K-manifold
\endsubhead
\subsubhead 3.1 Notation and preliminary facts
\endsubsubhead
\par
\bigskip
If $G$ is a compact semisimple Lie group and $\t^\C \subset \g^\C$ is 
a given Cartan subalgebra, we will  use the following notation:
\roster
\item"-" $\B$ is the Cartan-Killing form of $\g$ and for any subspace $A\subset \g$, 
$A^\perp$ is the $\B$-orthogonal complement to $A$;
\item"-" $R$ is the root system of $(\g^\C, \t^\C)$;
\item"-" $H_\alpha\in \t^\C$  is  the $\B$-dual element to the root $\alpha$; 
\item"-" for any $\alpha, \beta \in R$,  the scalar product $(\alpha, \beta)$ 
is set to be equal to $(\alpha, \beta) = \B(H_\alpha, H_\beta)$; 
\item"-" $E_\alpha$ is the root vector with root $\alpha$ in the Chevalley
normalization; in particular $\B(E_\alpha, E_{- \beta}) =  \delta_{\alpha\beta}$, 
$[E_\alpha, E_{-\alpha}] = H_\alpha$, $[H_\alpha, E_\beta] = (\beta, \alpha) E_\beta$ and 
$[H_{\alpha}, E_{-\beta}] = - (\beta, \alpha) E_{-\beta}$; 
\item"-" for any root $\alpha$, 
$$F_\alpha = \frac{1}{\sqrt{2}}(E_\alpha - E_{-\alpha})\ ,
\qquad G_\alpha = \frac{i}{\sqrt{2}}(E_\alpha + E_{-\alpha})\ ;$$
note that for $\alpha, \beta \in R$
$$\B(F_{\alpha}, F_{\beta}) = - \delta_{\alpha \beta} =  \B(G_{\alpha}, G_{\beta})\ ,\
 \B(F_{\alpha}, G_{\beta}) = \B(F_\alpha, H_\beta) = \B(G_\alpha, H_\beta) = 0\ ;$$
\item"-" the notation for the roots of a simple Lie algebra is the same of 
\cite{GOV} and \cite{AS}.
\endroster
Recall that for any two roots $\alpha, \beta$, with $\beta \neq - \alpha$, in case
 $[E_\alpha, E_\beta]$ is non trivial then it is equal to 
$[E_\alpha, E_\beta] = N_{\alpha, \beta} E_{\alpha + \beta}$ 
where the coefficients $N_{\alpha, \beta}$ verify the following 
conditions:
$$N_{\alpha, \beta} = - N_{\beta, \alpha}\ ,\quad  
N_{\alpha, \beta} = - N_{-\alpha, - \beta}\ .\tag 3.1$$
From (3.1) and the properties of root vectors in the Chevalley normalization, 
the following well known properties can be derived:
\roster
\item  for any  $\alpha, \beta \in R$ with $\alpha \neq \beta$
$$[F_\alpha, F_\beta], [G_\alpha, G_\beta] \in \span\{ F_\gamma\ ,\gamma \in R\}\ ,\qquad
[F_\alpha, G_\beta] \in \span\{ G_\gamma\ ,\gamma \in R\}\ ;\tag 3.2$$
\item for any $H \in \t^\C$ and any $\alpha, \beta \in R$,  $\B(H, [F_\alpha, F_\beta]) = 
\B(H, [G_\alpha, G_\beta]) = 0$ and 
$$\B(H, [F_\alpha, G_\beta]) = i \delta_{\alpha \beta} \B(H, H_\alpha) = \delta_{\alpha
\beta} \alpha(i H)\ ;\tag 3.3$$
\endroster
\bigskip
Finally, for what concerns the Lie algebra of flag manifolds and of 
CR manifolds, we adopt the following notation.\par
Assume that $G/K$ is a flag manifold with invariant complex 
structure $J$ (for definitions and basic facts,  we refer to
 \cite{Al}, \cite{AP}, \cite{BFR}, \cite{Ni}) and let
 $\pi: G/L \to G/K$ be a
$G$-equivariant  $S^1$-bundle over $G/K$. In particular, 
let us assume that $\l$
is a codimension one subalgebra of $\k$. Recall that $\k = \k^{ss} + \z(\k)$, with 
$\k^{ss}$ semisimple part of $\k$.  Hence  the semisimple part $\l^{ss}$ 
of $\l$ is equal to $\k^{ss}$ and  $\k = \l + \R Z = (\k^{ss} + \z(\k) \cap \l) + \R Z$
for some   $Z \in \z(\k)$. \par
Let  
$\t^\C\subset \k^\C$ be a Cartan subalgebra for $\g^\C$ contained in $\k^\C$ and 
$R$ the root system of $(\g^\C, \t^\C)$. Then we will use the following notation:
\roster
\item"-" $R_o = \{ \alpha \in R\ ,\ E_\alpha \in \k\ \}$;
\item"-" $R_\m = \{ \alpha \in R\ ,\ E_\alpha \in \m\ \}$;
\item"-" for any $\alpha \in R$, we denote by $\g(\alpha)^\C = \span_\C\{\ E_{\pm \alpha},
H_\alpha\}$ and $\g(\alpha) = \g(\alpha)^\C \cap \g$;
\item"-" $\m(\alpha)$ denotes the irreducible $\k^\C$-submoduli
of $\m^\C$, with  highest weight $\alpha\in R_\m$;
\item"-" if $\m(\alpha)$ and $\m(\beta)$ are equivalent
as $\l^\C$-moduli, we denote by $\m(\alpha) + \lambda\m(\beta)$  the irreducible 
$\l^\C$-module with  highest weight vector $E_\alpha + 
\lambda E_{\beta}$, $\alpha, \beta \in R_\m$, $\lambda\in \C$.
\endroster
\bigskip
\bigskip
\subsubhead 3.2 The structural decomposition $\g = \l + \R Z_\D + \m$ determined by the CR structure of
a regular
orbit
\endsubsubhead\par
\bigskip
The main results of this subsection are 
given by 
the following two theorems on the structural decomposition of the regular 
orbits of a K-manifolds. 
The first one is a straightforward consequence  of definitions, Theorem 2.6 and
the results
in \cite{PS}.\par
\medskip
\proclaim{Theorem 3.1} Let $(M, J, g)$ be a standard K-manifold acted on by the compact 
semisimple 
group $G$ and let 
$\g = \l + \R Z_\D + \m$ and $\m^{10}$ be the structural decomposition and the holomorphic 
subspace, respectively, 
associated with the CR structure $(\D, J)$ of a regular orbit $G/L = G\cdot p$. Let also 
$J: \m \to \m$ be the unique complex structure on $\m$, which determines the 
decomposition $\m^\C = \m^{10} + \overline{\m^{10}}$.\par
Then, $\k = \l + \R Z_\D$ is the isotropy subalgebra of a flag manifold $K$, 
and the complex structure 
$J: \m \to \m$ is $\operatorname{ad}_\k$-invariant and  corresponds to a
$G$-invariant complex structure $J$ on $G/K$.\par
In particular, 
there exists a Cartan subalgebra $\t^\C \subset \k^\C$ and an ordering of
the associated root system $R$, so that $\m^{10}$ is generated
by the corresponding positive root vectors in $\m^\C = (\k^\perp)^\C$.
\endproclaim
\bigskip
The following theorem describes the structural decomposition and the holomorphic subspace 
of a regular orbit of a non-standard K-manifold. Also this theorem can be considered as 
a consequence of Theorem 2.6, but the proof is  a little bit more involved.\par
\medskip
\proclaim{Theorem 3.2} Let $(M, J, g)$ be a non-standard K-manifold acted 
on by the compact semisimple 
 group $G$ and  let 
$\g = \l + \R Z_\D + \m$ and $\m^{10}$ be the structural decomposition and the holomorphic 
subspace, respectively, 
 associated with the CR structure $(\D, J)$ of a regular orbit $G/L = G\cdot p$. \par
Then there exists a simple  subalgebra $\g_F \subset \g$ with the following properties: 
\roster
\item"a)" denote by $\l_F = \l \cap \g_F$, $\l_o = \l \cap \g_F^\perp$, 
 $\m_F = \m \cap \g_F$ and $\m' = \m \cap \g_F^\perp$; then the pair $(\g_F, \l_F)$
is one of those listed in Table 1 and $\g$ and $\g_F$
admit  the following  $\B$-orthogonal decompositions:
$$\g = \l_o + (\l_F + \R Z_\D) + (\m_F + \m')\ ,\qquad \g_F = \l_F + \R Z_\D + \m_F\  ;$$
furthermore 
$[\l_o, \g_F] = \{0\}$ and the connected subgroup 
$K \subset G$ with Lie algebra $\k = \l_o + \g_F$ is the isotropy
subalgebra of a flag manifold
$G/K$; 
\item"b)"  denote by  $\m^{10}_F = \m_F^\C \cap \m^{10}$; then there
exists a   
Cartan subalgebra $\t^\C_F \subset \l^\C_F + \C Z_\D$ and a complex number $\lambda$
with $0< |\lambda| < 1$ so that  the element $Z_\D$, determined up to 
scaling, and the subspace $\m^{10}_F$, determined up to an element of the Weyl group 
 and  up to complex conjugation,  are as listed in Table 1
(see  \S 3.1 for  notation): 
\medskip
\moveleft 0.2cm
\vbox{\offinterlineskip
\halign {\strut\vrule\hfil\ $#$\ \hfil
 & \vrule\hfil\ $#$\
\hfil&\vrule\hfil\  $#$
\hfil &
\vrule
\hfil\  $#$\ 
\hfil\vrule \cr 
\noalign{\hrule} 
\phantom{\frac{\frac{1}{1}}{\frac{1}{1}}}
\g_F\ \ 
&
\l_F
& 
Z_\D
& \m^{10}_F
\cr 
\noalign{\hrule}
\goth{su}_{2}
& \{0\}
& 
-\frac{i}{2} H_{\varepsilon_1 - 
\varepsilon_2}
&
\underset{\phantom{A}}\to
{\overset{\phantom{B}}\to{
\C(E_{\varepsilon_1-\varepsilon_2} + \lambda
 E_{-\varepsilon_1+\varepsilon_2})
}}
\cr 
\noalign{\hrule}
\goth{su}_{n+1}
& \goth{su}_{n-2} \oplus \R
& 
-\ i H_{\varepsilon_1 - 
\varepsilon_2}
&
\underset{\phantom{A}}\to
{\overset{\phantom{B}}\to{\smallmatrix
(\C(E_{\varepsilon_1-\varepsilon_2} + \lambda^2
 E_{-\varepsilon_1+\varepsilon_2})\oplus \\
(\m(\varepsilon_1 - \varepsilon_3) + 
\lambda
 \m(\varepsilon_2 - \varepsilon_3)) \oplus
(\m(\varepsilon_3 - \varepsilon_2) + 
\lambda
 \m(\varepsilon_3 - \varepsilon_1)) 
\endsmallmatrix 
}}
\cr 
\noalign{\hrule}
\goth{su}_2 + \goth{su}_2
& \R 
& 
- \frac{i}{2}( H_{\varepsilon_1 - 
\varepsilon_2} + H_{\varepsilon'_1 - 
\varepsilon'_2})
&
\underset{\phantom{A}}\to
{\overset{\phantom{B}}\to{
\smallmatrix
\C(E_{\varepsilon_1 - 
\varepsilon_2} + \lambda E_{-(\varepsilon'_1 - 
\varepsilon'_2)}) 
\oplus \C(E_{\varepsilon'_1 - 
\varepsilon'_2} + \lambda E_{-(\varepsilon_1 - 
\varepsilon_2)}) 
\endsmallmatrix
}}
\cr \noalign{\hrule}
\goth{so}_7 
& \goth{su}_3
& 
- \frac{2 i}{3}( H_{\varepsilon_1 +
\varepsilon_2} + H_{\varepsilon_3})
&
\underset{\phantom{A}}\to
{\overset{\phantom{B}}\to{
\smallmatrix \left(
\m(\varepsilon_1 + \varepsilon_2) +
\lambda \m( - \varepsilon_3) \right) \oplus
\overline{\m( - \varepsilon_3)}
+  \left(\lambda \overline{\m(\varepsilon_1 + \varepsilon_2)}
 \right)
\endsmallmatrix
}}
\cr \noalign{\hrule}
\goth{f}_4 
& \goth{so}_7 
& 
- i 2 H_{\varepsilon_1}
&
\underset{\phantom{A}}\to
{\overset{\phantom{B}}\to{\smallmatrix
(\m(\varepsilon_1+\varepsilon_2) + \lambda^2
 \m(-\varepsilon_1+\varepsilon_2))\oplus \\
(\m(1/2(\varepsilon_1 +\varepsilon_2+ 
\varepsilon_3+\varepsilon_4)) +
\lambda
 \m(1/2(-\varepsilon_1+\varepsilon_2 + 
\varepsilon_3 +\varepsilon_4)))
\endsmallmatrix 
}}
\cr \noalign{\hrule}
\goth{so}_{2n+1} 
& \goth{so}_{2n-1} 
& 
- i H_{\varepsilon_1}
&
\underset{\phantom{A}}\to
{\overset{\phantom{B}}\to{\smallmatrix
(\m(\varepsilon_1+\varepsilon_2) + \lambda
 \m(-\varepsilon_1+\varepsilon_2))
\endsmallmatrix 
}}
\cr \noalign{\hrule}
\goth{so}_{2n} 
& \goth{so}_{2n-2}
& 
-i H_{\varepsilon_1}
&
\underset{\phantom{A}}\to
{\overset{\phantom{B}}\to{\smallmatrix
(\m(\varepsilon_1+\varepsilon_2) + \lambda
 \m(-\varepsilon_1+\varepsilon_2))
\endsmallmatrix 
}}
\cr \noalign{\hrule}
\goth{sp}_n 
& \smallmatrix
\goth{sp}_1 + 
\goth{sp}_{n-2} 
\endsmallmatrix
& 
- i H_{\varepsilon_1 + \varepsilon_2}
&
\underset{\phantom{A}}\to
{\overset{\phantom{B}}\to{\smallmatrix
(\m(2\varepsilon_1) + \lambda^2
 \m(-2\varepsilon_2))\oplus
(\m(\varepsilon_1+ \varepsilon_3) + \lambda
 \m(-\varepsilon_2 +\varepsilon_3))
\endsmallmatrix
}}
\cr \noalign{\hrule}
}}
\centerline{\bf Table 1}
\medskip
\item"c)" the holomorphic subspace $\m^{10}$ admits the following orthogonal
decomposition
$$\m^{10} = \m^{10}_F + \m'{}^{10}$$
where $\m'{}^{10} = \m'{}^\C \cap \m^{10}$; 
\item"d)"  the complex structure $J' : \m' \to \m'$ associated with the 
eigenspace decomposition 
$\m'{}^\C = \m'{}^{10} + \m'{}^{01}$, where $\m'{}^{01} = \overline{\m'{}^{10}}$,  is 
$\operatorname{Ad}_{K}$-invariant and 
determines a $G$-invariant complex structure on  the 
flag manifold $G/K$; in particular the $J'$-eigenspaces are $\operatorname{ad}_{\R Z_\D}$-invariant: 
$$[\R Z_\D, \m'{}^{10}] 
\subset \m'{}^{10}\ ,\qquad [\R Z_\D, \m'{}^{01}] \subset \m'{}^{01}\ .$$
\endroster 
\endproclaim
\bigskip
The proof of Theorem 3.2 needs the following  Lemma. \par
\bigskip
\proclaim{Lemma 3.3} Let  $G/L = G\cdot p$ be a regular 
orbit of the non-standard K-manifold $(M, J, g)$. 
Let also $\pi: (G/L, \D, J) \to (G/K, J_o)$ be the CRF fibration given in Theorem 2.6 
and   $(\D^K, J^K)$ the  CR structures of 
the standard fiber $K/L$. \par
Then: 
\roster
\item"i)" the 1-dimensional 
subspaces $\R Z_{\D^K}$ and $\R Z_\D$ of 
the  structural decompositions of $\k$ and 
$\g$ at the point $p$ 
are the same, i.e. their structural decompositions are
$\k = \l + \R Z_\D + \m_K$ and $\g = \l + \R Z_\D + \m = \l + \R Z_\D + (\m_K + \m')$; 
\item"ii)" the holomorphic subspace $\m^{10}$ of $(G/L, \D, J)$ admits the $\B$-orthogonal 
decomposition
$\m^{10} = \m^{10}_K + \m'{}^{10}$
where $\m'{}^{10} = \m^{10} \cap \m'{}^\C$ and $\m^{10}_K$ is the holomorphic subspace of $(K/L, \D^K, 
J^K$);
\item"iii)" $[\R Z_\D, \m'{}^{10}] \subset \m'{}^{10}$ and $[\R Z_\D, \m'{}^{01}] \subset \m'{}^{01}$.
\endroster
\endproclaim
\demo{Proof} Let $\k = \l + \R Z_{\D^K} + \m_K$ and $\g = \l + \R Z_\D + \m$
be the structural decompositions of $\k$ and $\g$ at the point $p$,
associated with the CR structures
$(\D^K, J^K)$ and $(\D, J)$, respectively. Denote also by $J^K$ and $J$ the 
induced complex structures on $\m_K$ and $\m$. \par
To prove i), we  have to show that 
 $\R Z_{\D^K} = \R Z_\D$. This  is proved by  the following observation. By definitions, 
$$\m_K= \{\ X\in \m\ :\ \pi_*(\hat X_{eL}) = 0\} = \m \cap \k$$
and hence  
$$\R Z_{\D^K} = \k \cap (\l + \m_K)^\perp = \k \cap (\l + (\m\cap \k))^\perp \subseteq 
\k \cap (\l + \m)^\perp = \k \cap \R Z_\D = \R Z_\D\ .$$
ii)  follows from   the fact  $J_K = J|_{\m_K}$.\par
To prove iii), we recall that  by Proposition 2.5, if $P$ is the parabolic subgroup 
such that $(G/K, J_F)$ is $G$-equivariantly biholomorphic to  $G^\C/P$, 
then the  subalgebra $\p = Lie(P) \subset \g^\C$
  verifies
$$\l^\C + \m^{01}_K + \m'{}^{01} \subset \p = \k^\C + \n$$
where $\n$ is the nilradical of $\p$ and $\k^\C$ is a  reductive  complement to $\n$. 
In particular, $\m'{}^{01} \subset \p \cap (\k^\C)^\perp = \n$. Moreover, 
$$\dim_\C \m'{}^{01} = 
\dim_\C G/P = \dim_\C \n$$ 
and hence $\m'{}^{01} = \n$. It follows that
$[\R Z_\D, \m'{}^{01}] \subset [\k^\C, \n] \subset \n = \m'{}^{01}$ and   
$[\R Z_\D, \m'{}^{10}] = \overline{[\R Z_\D, \m'{}^{01}]} \subset 
\overline{\m'{}^{01}} = \m'{}^{10}$.
\qed
\enddemo
\bigskip
\bigskip
\demo{Proof of Theorem 3.2} Let $K \subset G$
be a subgroup so that any regular orbit $G/L$ admits a very nice 
CRF fibration $\pi: (G/L, \D, J) \to (G/K, J_o)$
as prescribed by Theorem 2.6. Then, for any regular point $p$, 
the  $K$-orbit $K/L = K\cdot p \subset G/L = G \cdot p$ (which is the fiber
of the CRF fibration $\pi$) is a non-trivial Morimoto-Nagano space. 
In particular, $K/L$ is Levi non-degenerate,
it is simply connected and the CR structure is non-standard (for the definition of 
non-standard CR structures and the properties of 
the CR structures of the Morimoto-Nagano spaces,
see \cite{AS}). 
Furthermore, by Lemma 3.3, the 1-dimensional subspace $\R Z_{\D^K}$ associated with the 
 CR structure of $K/L$ coincides with the 1-dimensional subspace  $\R Z_\D$ associated
with the CR structure of $G/L$.\par
Let $L_o \subset L$ be the normal subgroup of the elements which act trivially on 
$K/L$. Let also  $G_F = K/L_o$ and 
$\l_o = Lie(L_o)$,  $\g_F = \k \cap (\l_o)^\perp \cong  Lie(G_F)$.\par
Note that Th. 1.3, 1.4 and 1.5 of \cite{AS} apply immediately to the homogeneous CR manifold
$G_F/L_F$, with $L_F = L\mod L_o$. In particular, since the CRF fibration 
$\pi: G/L \to G/K$ is very nice,  
$K/L = G_F/L_F$  is a
primitive homogeneous CR manifold (for the definition of primitive CR manifolds, see \cite{AS}) 
and $\g_F$  is  $\goth{su}_n$, 
$\goth{su}_2 + \goth{su}_2$, $\goth{so}_7$, $\goth{f_4}$, 
$\goth{so}_{n}$ ($n\geq 5$) or $\goth{sp}_n$ ($n\geq 2$).\par
From
Th.1.4, Prop. 6.3 and Prop. 6.4 in \cite{AS} and from Lemma 3.3 i) and ii), 
 it follows immediately that the subalgebra 
$\g_F$  and the holomorphic subspace $\m^{10}_F$, associated with 
the CR structure of the fiber $K/L  = G_F/L_F$, verify a), b), c) and d). \qed
\enddemo
\bigskip
In the following, we will call the subalgebra  $\g_F$ {\it  the Morimoto-Nagano
subalgebra of the non-standard K-manifold $M$\/}. We will soon prove that 
the Morimoto-Nagano
 subalgebra is independent (up to conjugation) from the choice of the 
regular orbit $G\cdot p = G/L$. \par
We will also 
call  $(\g_F,\l_F)$ and the subspace $\m^{10}_F$ the {\it Morimoto-Nagano 
pair\/} and the {\it Morimoto-Nagano holomorphic subspace\/}, respectively,  
of the regular orbit $G/L = G\cdot p$.\par
\bigskip 
\bigskip
\subsubhead 3.3 Optimal transversal curves 
\endsubsubhead
\bigskip
We prove now the existence of a special family of
 nice transversal curves called optimal transversal curves
(see \S 1). We first show the existence of such 
curves  for a non-standard K-manifold.\par
\bigskip
\proclaim{Theorem 3.4} Let $(M, J , g)$ be a non-standard K-manifold
acted on by the compact semisimple group $G$. Then there exists 
a point $p_o$ in the non-complex singular orbit and an element $Z \in \g$,
such that the curve 
$$\eta: \R \to M \ ,\qquad \eta_t = \exp(t i Z) \cdot p_o$$
verifies the following properties:
\roster
\item it is a nice transversal curve; in particular 
the isotropy subalgebra $\g_{\eta_t}$ for any 
$\eta_t\in M_{\text{reg}}$ is a fixed subalgebra $\l$;
\item there exists a subspace $\m$ such that, for any 
$\eta_t\in M_{\text{reg}}$, the structural decomposition 
$\g = \l + \R Z_\D(t) + \m(t)$ of the orbit $G/L = G\cdot \eta_t$ 
is given by  $\m(t) = \m$ and $\R Z_\D(t) = \R Z$;
\item the Morimoto-Nagano pairs $(\g_F(t), \l_F(t))$ of the regular orbits $G\cdot \eta_t$
do not dependent on $t$; 
\item for any $\eta_t \in M_{\text{reg}}$, the holomorphic subspace $\m^{10}(t)$ 
admits the orthogonal decomposition
$$\m^{10}(t) = \m^{10}_F(t) + \m'{}^{10}(t)$$
where $\m'{}^{10}(t) = \m'{}^{10} \subset \m^\C$  is independent on $t$ and 
 $\m^{10}_F(t)$ is a Morimoto-Nagano
holomorphic subspace which is listed  in Table 1, determined by 
the parameter $\lambda$ equal to
$$\lambda = \lambda(t) = e^{2t}\ .$$
\endroster
Moreover, if  $\eta_t =  \exp(t i Z) \cdot p_o$ is any of such curves and if $(\g_F, \l_F)$
is (up to conjugation) the Morimoto-Nagano pair of a regular orbits $G/L = G\cdot \eta_t$, 
then (up to conjugation) $Z$ is  the element in the column ''$Z_\D$''
of Table 1, associated with the Lie algebra 
$\g_F$. 
\endproclaim
\bigskip
For the proof of Theorem 3.4, we first need two Lemmata. \par
\bigskip
\proclaim{Lemma 3.5} Let $(M, J, g)$ be a K-manifold acted on 
by the compact semisimple Lie group $G$. Let also  $p$ be a regular point and 
$G/L = G\cdot p$ and $G^\C/ H = G^\C\cdot p$ the $G$- and the $G^\C$-orbit
of $p$, respectively. Then:
\roster
\item the isotropy subalgebra $\h = Lie(G^\C_p)$ is equal to 
$$\h = \l^\C + \m^{01}$$
where $\m^{01} = \overline{\m^{10}}$ is the anti-holomorphic subspace
associated with the CR structure  of $G/L = G\cdot p$;
\item   for any $g\in G^\C$,  the isotropy 
subalgebra $\l' = \g_{p'}$ at 
$p' = g\cdot p$ is equal to 
$$\l' = \operatorname{Ad}_g(\l^\C + \m^{01})\cap \g\ ;$$
\item let $g\in G^\C$ and suppose that  $p' = g \cdot p$ is a regular point; 
if we denote by $\g = \l' + \R Z'_\D + \m'$ and by  $\m'{}^{10}$
  the structural decomposition and the holomorphic subspace, respectively, 
given by the CR structure of  $G\cdot p' = G/L'$, then
$$
\m'{}^{10} = \overline{\operatorname{Ad}_g(\l^\C + \overline{\m^{10}})}\ ,$$
$$
\m' = \left(\operatorname{Ad}_g(\l^\C + \overline{\m^{10}}) + 
\overline{\operatorname{Ad}_g(\l^\C + \overline{\m^{10}})}\right)\cap \g \cap \l'{}^\perp\ .$$
\endroster
\endproclaim
\demo{Proof} (1) Consider an element $V = X + i Y\in \g^\C$, with $X, Y \in \g$. 
Then $V$  belongs to $\h$
if and only if 
$$\widehat{X + i Y}|_p = \hat X_p  + J \hat Y_p = 0\ .$$
This means that $J \hat X_p = - \hat Y_p$ is tangent to the orbit $G\cdot p$. In particular,
$X, Y \in \l + \m$ and    $V = X + i J X \in \l^\C + \m^{01}$.\par
(2) Clearly, $L' = G \cap G^\C_{p'} = G \cap (g H g^{-1})$ and 
$\l' = \g \cap \operatorname{Ad}_g(\h)$. The claim is then an immediate consequence of (1).\par
(3) From (1), it follows that 
$$\m'{}^{10} = \overline{\m'{}^{01}} = \overline{\h' \cap (\l'{}^\C)^\perp} = 
\overline{\operatorname{Ad}_g(\l^\C + \m^{01})}\cap (\l'{}^\C)^\perp\ .$$
From this, the conclusion follows.\qed
\enddemo
\bigskip
\proclaim{Lemma 3.6} Let $(M, J, g)$ be a  K-manifold acted on 
by the compact semisimple Lie group $G$. Let also $p$ be a regular point and 
$\g = \l + \R Z_\D + \m$ the structural decomposition  associated with the CR structure
of $G/L = G\cdot p$. Then:
\roster
\item for any $g \in \exp(\C^* Z_\D)$,    the isotropy subalgebra 
$\g_{p'}$ at the point 
$p' = g\cdot p$ is orthogonal to $\R Z_\D$; moreover, $\l \subseteq \g_{p'}$
and, if $p'$ is regular, 
$\l = \g_{p'}$; 
\item the curve 
$$\eta : \R \to M\ ,\qquad \eta_t = \exp(i t Z_\D)\cdot p$$
is a nice transversal curve through $p$.
\endroster
\endproclaim
\demo{Proof} (1) From Lemma 3.5 (2), for any point $p' = \exp(\lambda Z_\D)\cdot p$, with $\lambda \in \C^*$,
$$\B(\g_{p'},\R Z_\D) = \B(\operatorname{Ad}_{\exp(\lambda Z_\D)}(\l^\C + \m^{01})\cap \g, \R Z_\D) = $$
$$ =
\B((\l^\C + \m^{01})\cap \g, \operatorname{Ad}_{\exp(-\lambda Z_\D)}(\R Z_\D)) 
= \B((\l^\C + \m^{01})\cap \g,  \R Z_\D) = 0\ .$$
Moreover, since $Z_\D \in C_{\g^\C}(\l^\C)$, we get that 
$$\g_{p'} = (\operatorname{Ad}_{\exp(\lambda Z_\D)}(\l^\C + \m^{01})) \cap \g = \l + 
\operatorname{Ad}_{\exp(\lambda Z_\D)}(\m^{01}) \cap \g \supset \l\ .$$
This implies  that $\l = \g_{p'}$ if 
$p'$ is regular.\par
\medskip 
(2) From (1), we have that 
condition (2.1) and the equality 
$G\cdot {\eta_t} = G\cdot p = G/L$
are verified for any point $\eta_t \in M_{\text{reg}}$. 
 It remains to show that $\eta$ intersects any 
regular orbit. \par
Let $\Omega = M\setminus G$ be the
orbit space  and $\pi: M \to \Omega = M\setminus G$ the natural projection map.
It is known (see e.g. \cite{Br}) that $\Omega$ is 
 homeomorphic to  $\Omega = [0,1]$, with $M_{\text{reg}} = \pi^{-1}(]0,1[)$.
Hence $\eta$ intersects any regular orbit if and only if $(\pi\circ \eta)(\R) \supset ]0,1[$.\par
Let $x_1 = \inf (\pi\circ\eta)(\R) $ and let $\{t_n\} \subset ]0, 1[$ be  a sequence
such that  $(\pi\circ \eta)_{t_n}$ tends to $x_1$. If we assume that 
$x_1 > 0$,  we may select a subsequence  $t_{n_k}$ so that
$\lim_{n_k\to \infty} \eta_{t_{n_k}}$ exists and it is equal to 
 a regular point $p_o$.  From (1)
and a  continuity argument, we could conclude that 
$\l$ is equal to the isotropy subalgebra $\g_{p_o}$, that
$\hat Z_\D|_{p_o} \neq 0$ and that $J \hat Z_\D|_{p_o}$ is not tangent
to the orbit $G\cdot p_o$.  In particular, it would follow
 that the curve $\exp(i \R Z_\D) \cdot p_o$   
has non-empty intersection with $\eta(\R) =\exp(i \R Z_\D) \cdot p$ and  that
$p_o \in \eta(\R)$; moreover we would have  that $\eta$ is transversal to $G\cdot p_o$ and 
that   $x_1 = \pi(p_o)$ is an inner point of $\pi\circ \eta(\R)$, 
which is a contradiction.\par
A similar contradiction arises if we assume that  $x_2 = \sup \pi\circ\eta(\R)  < 1$.
\qed \enddemo
\bigskip
\demo{Proof of Theorem 3.4} Pick a regular point $p$. Let 
$\g = \l + \R Z_\D + \m$ be
the structural decomposition of the orbit $G\cdot p$  and let 
$\eta_t = \exp(i t Z_\D)\cdot p$. From Lemmata 3.5 and 3.6 and Theorem 3.2, 
 the structural decompositions $\g = \l + \R Z_\D(t) + \m(t)$ 
of all regular orbits $G\cdot \eta_t$ are independent on $t$.
Moreover, from Lemma 3.5  and  Theorem 3.2, it follows 
that the Morimoto-Nagano pair $(\g_F, \l_F)$ is the same for all 
regular orbits $G\cdot \eta_t$ and
the holomorphic subspace $\m^{10}_t$
 of the orbit $G\cdot \eta_t$ is of the 
form
$$\m^{10}_t = \overline{\operatorname{Ad}_{\exp(i t Z_\D)}(\overline{\m^{10}_0)}}  = 
\operatorname{Ad}_{\exp( - i t Z_\D)}(\m^{10}_F(0)) + 
\operatorname{Ad}_{\exp(- i t Z_\D)}(\m'{}^{10}(0))
\tag 3.4$$
where $\m^{10}_0 = \m^{10}_F{}(0) + \m'{}^{10}(0)
$ is the decomposition of the holomorphic 
subspace of $G\cdot \eta_0$ given in 
Theorem 3.2 c). Since $Z_\D\in \g_F$, 
from (3.4) and Theorem 3.2 d), it follows 
that
$$\m^{10}_t = \operatorname{Ad}_{\exp(-i t Z_\D)}(\m^{10}_F(0)) + \m'{}^{10}(0)\ .$$
This proves that the Morimoto-Nagano holomorphic subspace $\m^{10}_F(t)$ of the 
orbit $G\cdot \eta_t$ is 
$$\m^{10}_F(t) = \operatorname{Ad}_{\exp(- i t Z_\D)}(\m^{10}_F(0))\tag 3.5$$
and that the $\B$-orthogonal complement $\m'{}^{10} =  \m'{}^{10}(0)$ is 
independent on $t$ and $\operatorname{ad}_{Z_\D}$-invariant.\par
A simple computation shows  that  if $\g_F$ and  $\m^{10}_F(t) = 
\operatorname{Ad}_{\exp(-i t Z_\D)}(\m^{10}_F(0))$ 
appear in a row of
 Table 1 and if $Z_\D$ is equal to $Z_\D= A Z_o$, where 
$Z_o$ is  the corresponding element 
listed  in the column ''$Z_\D$'', then 
$\m^{10}_F(t)$ is determined by a complex
parameter  $\lambda = \lambda(t)$, 
which verifies  the differential equation
$$\frac{d \lambda}{dt}  = 2 A\lambda(t)\ .$$
In particular, if we assume $A = 1$, then $\lambda(t) =  e^{2t + B_p}$ where
$B_p$ is a complex number which depends only on
 the regular point $p$. \par
Let us replace $p$ with the point $p_o = \exp( - i \frac{B_p}{2} Z) \cdot p$: it is immediate
to realize that  the new function $\lambda(t)$ is  equal to
$$\lambda(t) =  e^{2t + B_p - B_p} = e^{2t}\ .$$
This proves that the curve $\eta_t = e^{i t Z_\D}\cdot p_o$ verifies (1), (2), (3) and (4). \par
It remains to prove that for any choice of the regular
point $p$,  the point  $p_o = \exp( - i\frac{B_p}{2} Z)\cdot  p$ is  a 
point of the non-complex singular orbit of $M$. \par
Observe that, since $\eta(\R)$
is the orbit of a real 1-parameter subgroup of $G^\C$, the complex isotropy  subalgebra 
$\h_t \subset \g^\C$ is (up to conjugation) independent on the point $\eta_t$. Indeed, if 
$\eta_{t_o}$ is a regular point 
with complex isotropy subalgebra $\h_{t_o} = \l^\C + \m^{01}_F + \m'{}^{01}$, 
then for any other point $\eta_t$, we have that
$$\h_{t} = \operatorname{Ad}_{\exp(i (t - t_o) Z_\D)}(\l^\C + \m^{01}_F + \m'{}^{01})\ .$$
On the other hand, the real isotropy subalgebra   
$\g_{\eta_t} \subset \g$ is equal to 
$$\g_{\eta_t} = \h_t \cap \g = 
\operatorname{Ad}_{\exp(i (t - t_o) Z_\D)}(\l^\C + 
\m^{01}_F + \m'{}^{01})\cap \g\ .\tag 3.6$$
From (3.6), Table 1 and (4), 
one can  check that in all cases
$$\g_{\eta_0} \supsetneq \l + \R Z_\D$$
and hence that $\eta_0 = p_o$ is a singular point for the $G$-action. 
On the other hand $p_o$ cannot 
be in the complex singular $G$-orbit, because otherwise this orbit would 
coincide with  $G^\C \cdot p_o = G^\C\cdot p$ and it
would contradict the  assumption 
that $p$  is a regular point for the 
$G$-action.
\qed
\enddemo
\bigskip
The following
is the analogous 
result for standard K-manifolds. \par
\bigskip
\proclaim{Theorem 3.7} Let $(M, J , g)$ be a standard K-manifold
acted on by the compact semisimple group $G$ and let $p_o$
be any regular point for the $G$-action. Let also $\g = \l + \R Z  + \m$
and $\m^{10}$  be the structural decomposition and the holomorphic 
subspace associated with the CR structure of the orbit $G/L = G\cdot p_o$. Then the curve 
$$\eta: \R \to M \ ,\qquad \eta_t = \exp(t i Z) \cdot p_o$$
verifies the following properties:
\roster
\item it is a nice transversal curve; in particular the stabilizer in $\g$ of any regular point 
$\eta_t$ is equal to the isotropy  subalgebra $\l = \g_{p_o}$;
\item  for any regular point $\eta_t$,  the structural decomposition 
$\g = \l + \R Z_\D(t) + \m(t)$ and the holomorphic subspace $\m^{10}(t)$
of the CR structure of $G/L = G\cdot \eta_t$ 
is given by the subspaces $\m(t) = \m$, $\R Z_\D(t) = \R Z$ and $\m^{10}(t) = \m^{10}$.
\endroster
\endproclaim
\demo{Proof} (1) is immediate from Lemma 3.6.\par
(2) It  is sufficient to prove  that $[Z, \m^{10} ] \subset \m^{10}$. In fact, 
from this
the claim follows as an immediate corollary of Lemmata 3.5 and 3.6. \par
Let $(G/K, J_F)$ be the flag manifold with invariant 
complex structure $J_F$, given by Theorem 2.6, so that  any 
regular orbit $G\cdot x$ admits a CRF fibration onto $G/K$, with fiber $S^1$.
 Let also 
$P$ be the parabolic subalgebra of $G^\C$ such that $G/K$ is biholomorphic to
$G^\C/P$. \par
From Proposition 2.5, if we denote by $\p = \k^\C + \n$ the decomposition 
of the parabolic subalgebra 
$\p \subset \g^\C$ into nilradical $\n$ plus reductive part $\k^\C$, we have that 
$$\k = \p \cap \g\ ,\qquad \l^\C \subsetneq \k^\C \ , 
\qquad \l^\C + \m^{01} \subset \k^\C + \n\ .\tag 3.7$$
Since the CRF fibration has fiber $S^1$, it follows that 
$\k = \l + \R Z'$
for some $Z' \in \z(\k) \subset \a = 
C_\g(\l) \cap \l^\perp$.  \par
In  case $\dim\a = 1$,  we have that $\a = \R Z  = \R Z'$ and hence 
 $\m^{10} \subset (\l^\C + \C Z)^\perp = (\k^\C)^\perp$.  
From (3.7)  we get  that $\m^{01} = \n$  and  that
$[Z, \m^{01}] \subset [\k^\C, \n] \subset \n = \m^{01}$.\par
In case $\a$ is 3-dimensional, 
let us denote by $\a^\perp = \a \cap \m = \a \cap (\R Z)^\perp$ and 
by $\a^{10} = \a^\C \cap \m^{10}$, $\a^{01} = \a^\C \cap \m^{01}
= \overline{\a^{10}}$ so that $(\a^\perp)^\C = 
\a^{10} + \a^{01}$. Consider also the orthogonal decompositions
$$\g = \l + \R Z + \m = \l + \R Z + \a^\perp + \m' \ ,
\qquad \m^{10} = \a^{10} + \m'{}^{10}\ ,$$
where $\m'{}^{10} = \m^{10} \cap \m'{}^\C$. 
Let $\l^{ss}$ be the semisimple part of $\l$ and note that 
$\l^{ss} = \k^{ss}$. 
By classical properties of flag manifolds (see e.g. \cite{Al},
\cite{AP}, \cite{Ni}) the $\operatorname{ad}_{\k^{ss}}$-module
 $\m'$ contains no trivial $\operatorname{ad}_{\k^{ss}}$-module and hence
$\m'{}^{10} = [\k^{ss}, \m'{}^{10}] = [\k, \m'{}^{10}]$. 
In particular, $\m'{}^{01} = \overline{\m'{}^{10}}$
is orthogonal to $\k^\C$ and hence it is included in $\n$. So,
$$[Z, \m'{}^{01}] \subset [Z, \n \cap (\l^\C + \a^\C)^\perp] 
\subset \n \cap (\l^\C + \a^\C)^\perp = 
\m'{}^{01}\ .$$
From this, it follows that 
 in order to prove that $[Z, \m^{10} ] \subset \m^{10}$, 
one has only to show that 
$[Z, \a^{10} ] \subset \a^{10} \subset \m^{10}$. \par
By dimension counting, $\a^{10} = \C E$
for some element $E\in \a^\C \simeq \goth{sl}_2(\C)$. In case $E$ is a nilpotent element for 
the Lie algebra $\a^\C \simeq \goth{sl}_2(\C)$, we may choose a Cartan subalgebra $\C H_\alpha$
for $\a = \goth{sl}_2(\R)$, 
so that $E \in \C E_{\alpha}$. In this case, we have that
 $$Z \in (\a^{10} + \a^{01})^\perp = 
(\C E_\alpha + \C E_{-\alpha})^\perp = \C H_\alpha$$
 and hence
$[Z, \a^{10}] \subset [\C H_\alpha, \C E_\alpha ] = \C E_\alpha = \a^{10}$ and we are done.\par
In case $E$ is a regular element for $\a^\C$,  with no loss of generality, 
we may consider a Cartan subalgebra $\C H_\alpha$ for $\a^\C$ 
so that $\C E = \C (E_\alpha + t E_{-\alpha})$
for some $t\neq 0$. In this case, $\a^{01} = \overline{\a^{10}} = \C (E_{-\alpha} + \bar t E_{\alpha}) = 
\C(E_{\alpha} + \frac{1}{\bar t} E_{-\alpha})$ 
and, since $\a^{10} \cap \a^{01} = \{0\}$, it follows that $t \neq 1/\bar t$. 
In particular, we get that
$\C Z = (\a^{10} + \a^{01})^\perp = \C H_\alpha$. Now, by 
Lemma 3.5 (2), for any $\lambda \in \C^*$, the isotropy subalgebra $\l_{g_\lambda\cdot p_o}$,
with $g_{\lambda} = \exp(\lambda Z) $, is equal to 
$$\l_{g_{\lambda}\cdot p_o} = Ad_{\exp(\lambda Z)} (\l^\C + \a^{01} + \m'{}^{01})\cap \g = 
\l^\C + \m'{}^{01} + \C (E_\alpha + t e^{- 2 \lambda \alpha(Z)} E_{-\alpha}) \cap \g\ .$$
Therefore, if $\lambda$ is such that $t e^{- 2 \lambda \alpha(Z)} = - 1$, we have that 
$\l_{g_{\lambda}\cdot p_o} = \l + \R (E_\alpha - E_{-\alpha}) \supsetneq \l$ and hence that 
$p = g_{\lambda} \cdot p_o$ is a singular point for the $G$-action. On the other hand, $p$ is in 
the $G^\C$-orbit of $p_o$ and hence the singular orbit $G\cdot p$ is not a complex
orbit. But this is in contradiction with the hypothesis that $M$ is standard and hence that 
it has  two  singular $G$-orbits, which are both complex.
\qed
\enddemo
\bigskip
Any curve $\eta_t = \exp(i t Z)\cdot p_o$, which verifies the 
claim of Theorems 3.4 or 3.7,  
will be called {\it optimal transversal curve\/}. \par
\bigskip
\bigskip
\subsubhead 3.4 The optimal bases  along the optimal transversal curves
\endsubsubhead\par
\bigskip
In all the following,   $\eta$ is an optimal transversal curve. 
In case $M$ is a non-standard 
K-manifold, we   denote by $\g = \l + \R Z_\D + \m$, $(\g_F, \l_F)$, $\m^{10}_F(t)$
 and $\m^{10} = \m^{10}_F(t) + \m'{}^{10}$
the structural decomposition, the Morimoto-Nagano pair, 
the Morimoto-Nagano subspace and the holomorphic 
subspace, respectively, at the regular points $\eta_t \in M_{\text{reg}}$. The same notation 
will be adopted   in case $M$ is a 
standard K-manifold, with the convention that, in this case, the Morimoto-nagano pair
 $(\g_F, \l_F)$ is the trivial pair $(\{0\}, \{0\})$ and that the Morimoto-Nagano
holomorphic subspace is
$\m^{10}_F = \{0\}$.\par
\medskip
We will also assume that 
$\l = \l_o + \l_F$, where 
$\l_o = \l \cap \l_F^\perp$. By 
 $\t^\C = \t^\C_o + \t^\C_F \subset \l^\C\subset \g^\C$, with $\t_o \subset \l_o$ and
$\t_F \subset \l_F$, we denote a Cartan subalgebra of $\g^\C$ 
with the property that, the expressions of 
$\m^{10}_F(t)$ and $Z_\D $ in terms of the root vectors of $(\g^\C_F, \t^\C_F)$
are exactly  as those 
listed in Table 1, corresponding 
to the parameter $\lambda_t = e^{2t}$. 
\par
Let $R$ be 
the root system  of $(\g^\C, \t^\C)$. Then $R$ is union of 
the following disjoint subsets of roots:
$$R = R^o \cup R' = (R^o_\perp \cup R^o_F) \cup (R_F' \cup R'_+ \cup R'_-)\ ,$$
where 
$$R^o_\perp = \{\ \alpha \ ,\ E_{\alpha} \in \l^\C_o\ \}\quad ,\quad
R^o_F = \{\ \alpha \ ,\ E_{\alpha} \in \l^\C_F\ \}\ ,$$
$$R'_F = \{\ \alpha \ ,\ E_{\alpha} \in \m^\C_F\ \}\quad ,\quad
R'_+ = \{\ \alpha \ ,\ E_{ \alpha} \in \m'{}^{10}\ \}\quad, \quad 
R'_- = \{\ \alpha \ ,\ E_{ \alpha} \in \m'{}^{01}\ \}\ .$$
Note that 
$$ - R^o_\perp = R^o_\perp\ ,\quad - R^o_F = R^o_F \ , \quad - R'_F = R'_F \ ,\quad - R'_+ = R'_-\ .$$
Moreover,  $R^o_\perp$ is orthogonal to $R^o_F$ and  
$R^o_\perp$, $R^o_F$ and $R^o_F \cup R'_F$ are closed subsystems. \par
Clearly, in case $M$ is standard, we will assume that $R^o_F = R'_F = \emptyset$.\par
\medskip
We claim that for any $\alpha \in R'_F$ there exists exactly one root 
$\alpha^d \in R'_F$ and two integers 
$\epsilon_\alpha = \pm 1$ and $\ell_\alpha = \pm 1, \pm 2$ such that, 
for any $t\in \R$, 
$$E_\alpha + e^{2\ell_\alpha t} \epsilon_{\alpha}E_{ -\alpha^d} \in \m^{10}_F(t)\ ,
\ .\tag 3.8$$
The proof of this claim is the following. 
By direct inspection of Table 1,
the reader can check that 
any maximal $\l^\C_F$-isotopic subspace of $\m^\C_F(t)$ (i.e. any maximal subspace which 
is sum of equivalent irreducible $\l^\C_F$-moduli) is direct sum of exactly 
two irreducible $\l^\C_F$-moduli
(see also \cite{AS}). Let us denote by $(\alpha_i, -\alpha^d_i)$ ($i = 1, 2, \dots$)
all  pairs of roots in $R_F$ with the property that the associated root vectors
$E_{\alpha_i}$ and $E_{-\alpha^d_i}$
are maximal weight vectors of  equivalent $\l^\C_F$-moduli in $\m^\C_F(t)$. Using  Table 1, 
one can check that in all cases $\m^{10}_F(t)$  decomposes into non-equivalent 
irreducible $\l^\C_F$-moduli, with maximal weight vectors  of the form 
 $$E_{\alpha_i} + \lambda^{(i)}_t E_{-\alpha^d_i}$$
where $\lambda^{(i)}_t = (\lambda(t))^{\ell_i} = 
e^{t \ell_i t}$, where $\ell_i$ is an integer which  is  either $\pm1$ or $\pm2$. \par
Hence  $\m^{10}_F(t)$
 is spanned by the vectors 
$E_{\alpha_i} + \lambda_t E_{-\alpha^d_i}$ and by  vectors of the form
$$[E_\beta, E_{\alpha_i} + \lambda_t E_{-\alpha^d_i}] = 
N_{\beta, \alpha_i} E_{\alpha_i + \beta} + \lambda_t
N_{\beta, - \alpha^d_i} E_{-\alpha^d_i + \beta}\ ,\tag 3.9$$
for some $E_\beta \in \l^\C$. Since the $\l^\C$-moduli containing $E_{\alpha_i}$
and $E_{-\alpha^d_i}$ are equivalent, the lengths of the sequences of roots 
$\alpha_i + r \beta$ and $-\alpha_i^d + r\beta$ are both 
equal to some 
given integer, say   $p$. This implies that 
for any root $\beta \in R^o_F$
$$N_{\beta, \alpha_i}^2 = (p+1)^2 = N_{\beta, -\alpha_i^d}^2$$
and hence that $\frac{ N_{\beta, \alpha_i}}{N_{\beta, - \alpha^d_i}} = \pm 1$. 
From this remark and (3.9), we conclude that
 $\m^{10}_F(t)$ is generated by elements of the form 
$$E_{\alpha} + \epsilon_{\alpha} e^{t \ell_\alpha t}
 E_{-\alpha^d}\ ,$$
where $\beta \in R_F^o$, $\alpha = \alpha_i + \beta$, 
$\alpha = \alpha_i + \beta$, $\alpha^d = \alpha^d_i + \beta$ and 
$\epsilon_{\alpha} = \frac{ N_{\beta, \alpha_i}}{N_{\beta, - \alpha^d_i}}$.  
This concludes the proof of the claim.\par
\bigskip
For any root $\alpha \in R_F$, we call {\it CR-dual root of $\alpha$\/} the root
$\alpha^d$ so that $E_{\alpha} + \epsilon_{\alpha} e^{t \ell_\alpha t}
 E_{-\alpha^d} \in \m^{10}(t)$. 
\medskip
We  fix  a positive root subsystem $ R^+ \subset R$ so that 
$R'_+ = R^+ \cap (R \setminus (R^o \cup R^o_F \cup R'_F))$. Moreover, we decompose 
the set of roots $R'_F$ into 
$$R'_F = R^{(+)}_F \cup R^{(-)}_F$$ 
where 
$$R^{(+)}_F  = \{ \alpha\in \ R'_F\ : \ E_{\alpha} + 
\epsilon_\alpha e^{\ell_\alpha t} E_{- \alpha^d}\in \m^{10}\ , \ \text{with}\ \ell_\alpha = + 1, +2\ \}$$
$$R^{(-)}_F  = \{ \alpha\in \ R'_F\ : \ E_{\alpha} + 
\epsilon_\alpha e^{\ell_\alpha t} E_{- \alpha^d}\in \m^{10}\ , \ \text{with}\ \ell_\alpha = - 1, -2\ \}$$
Using Table 1, one can check that in all cases
$$\m^{10} = \text{span}_\C\{\ E_\alpha + \epsilon_\alpha e^{\ell_\alpha t} E_{- \alpha^d}\ ,\ 
\alpha \in R^{(+)}_F\ \}$$
and that if $\alpha \in R^{(+)}_F$, then also  the CR dual root $\alpha^d \in R^{(+)}_F$.
We will denote by
$\{\alpha_1,  \alpha^d_1,  \alpha_2, \alpha^d_2, \dots , \alpha_r, \alpha^d_r\}$
the set of roots in $R^{(+)}_F $ and by $\{ \beta_1, \dots, \beta_s\}$
 the  roots in $R'_+ = R^+\cap R'$. \par
\medskip
Observe that the number  of  roots in $R^{(+)}_F$ is equal to
$\frac{1}{2}\left(\dim_\R G_F/L_F - 1\right)$, 
where $G_F/L_F$ is the Morimoto-Nagano space 
associated with the pair $(\g_F, \l_F)$.\par
\medskip
Finally, we consider 
the following basis for $\R Z_\D + \m \simeq T_{\eta_t} G\cdot \eta_t$.
We set
$$F_0 = Z_\D\ ,$$
and, for any $1 \leq  i \leq  r$, 
 we define the vectors $F^+_{i}$, $F^-_{i}$, $G^+_i$ and $G^-_i$, as follows: in case
$\{\alpha_i, \alpha_i^d\}\subset R^{(+)}_F$ is  a pair of CR dual roots with $\alpha_i \neq \alpha_i^d$,
we set 
$$F^+_{i} = \frac{1}{\sqrt{2}}(F_{\alpha_i} + \epsilon_{\alpha_i} F_{\alpha^d_i})\ ,
\quad F^-_{i} = \frac{1}{\sqrt{2}}(F_{\alpha_i} - \epsilon_{\alpha_i} F_{\alpha^d_i})\ ,$$
$$G^+_{i} = \frac{1}{\sqrt{2}}(G_{\alpha_i} + \epsilon_{\alpha_i} G_{\alpha^d_i})\ ,
\quad G^-_{i} = \frac{1}{\sqrt{2}}(G_{\alpha_i} - \epsilon_{\alpha_i} G_{\alpha^d_i})\ ,
\tag3.10$$
where $\epsilon_{\alpha_i} = \pm 1$ is the  integer which is defined in (3.8); in case
$\{\alpha_i, \alpha_i^d\}\subset R^{(+)}_F$ is  a pair of CR dual roots with $\alpha_i  =  \alpha_i^d$, 
we set 
$$F^+_{i} = F_{\alpha_i} = \frac{E_{\alpha_i} - E_{-\alpha_i}}{\sqrt{2}}\ ,\quad
G^+_{i} = G_{\alpha_i} = i\frac{E_{\alpha_i} + E_{-\alpha_i}}{\sqrt{2}}\tag 3.10'$$
and   {\it we  do not define the corresponding 
vectors $F^-_{i}$ or $G^-_{i}$\/}. Finally, for any $1\leq i\leq s = n-1 - 2 r$, we set
$$F'_{i} = F_{\beta_i}\ ,\quad
G'_{i} = G_{\beta_i}\ .
\tag3.11$$
Note that in case  $r$ is odd,  there is only one root $\alpha_i \in R^{(+)}_F$ such  
that $\alpha_i = \alpha^d_i$. When $\g_F = \su_2$, this root is also  the 
{\it unique\/} root in $R^{(+)}_F$.
\par
In case $\g_F = \{0\}$, we  set $F_0 = Z_\D$ and $F'_{i} = F_{\beta_i}$, 
$G'_{i} = G_{\beta_i}$ and we do not define the vector $F^{(\pm)}_i$ or 
$G^{(\pm)}_i$.\par
\medskip
The basis $(F_0, F^{\pm}_k, F_j, G^{\pm}_k, G_j)$  for $\R Z_\D + \m$, which we just defined, 
  will be called 
{\it optimal  basis  associated with the optimal  transversal 
curve $\eta$\/}. Notice that  this basis is $\B$-orthonormal.\par
\medskip
 For simplicity of notation, we will often use the symbol $F_k$ (resp. $G_k$) 
to denote any vector in the set  $\{F_0, F^{\pm}_j, F'_j\}$ (resp. 
in $\{G^{\pm}_j, G'_j\}$). 
We will also denote by  $N_F$ the number of elements of the form $F^{\pm}_i$. Note that 
$N_F$  is equal to half the real  dimension of the holomorphic
distribution of the Morimoto-Nagano space $G_F/L_F$. \par
For any  odd integer  $1 \leq 2k - 1 \leq N_F$, 
 we will assume that $F_{2k - 1} = F^+_k$; for any even 
integer   $2 \leq 2k \leq N_F$, we will assume $F_{2k} = F^-_k$. If $N_F$ is odd, 
we denote by $F_{N_F}$ the unique vector defined by (3.10').
We will also assume that $F_j = F'_{j - N_F}$ for any $N_F+1 \leq j\leq n-1$. \par
In case  $M$ is a standard K-manifold, we assume that 
 $N_F = 0$.\par
\bigskip
In the following lemma, we describe the action of the complex structure $J_t$
in terms of an optimal basis.\par
\medskip
\proclaim{Lemma 3.8} Assume that $\eta_t$ is an optimal transversal curve and
let 
$$(F_0, F^{\pm}_k, F'_j, G^{\pm}_k, G'_j)$$
  an associated optimal basis of 
$\R Z_\D + \m$. Let also  $J_t$ be the complex structure of $\m$
corresponding to the CR structure of a regular orbit $G\cdot \eta_t$. \par
Then $J_t F'_i = G'_i$ for any  $1 \leq i \leq s = n-1 - N_F$. Furthermore, if 
 $M$ is  non-standard (i.e. $N_F >0$) then:
\roster
\item if $1 \leq i\leq N_F$ and $\{\alpha_i, \alpha^d_i\}$,   is a pair of 
CR-dual roots in $R^{(+)}_F$ with $\alpha_i \neq \alpha^d_i$ then 
$$J_t F^+_i = - \coth(\ell_i t) G^+_i\ ,\quad J_t F^-_i = - \tanh(\ell_i t) G^-_i\ ,
\tag 3.12$$
where $\ell_i$ is equal to 
$2$ if $F^\pm_i \in [\m_F, \m_F]^\C \cap \m_F^\C$ and 
it is equal to $1$ otherwise; 
\item if $1 \leq i\leq N_F$ and  $\{\alpha_i, \alpha^d_i\}$  is a pair of 
CR-dual roots in $R^{(+)}_F$ with $\alpha_i =\alpha^d_i$, so that $F^+_i = F_{\alpha_i}$,
 then 
$$J_t F^+_i = - \coth(\ell_i t) G^+_i\ ,\tag 3.13$$ 
where $\ell_i$ is equal to 
$2$ if $F^\pm_i \in [\m_F, \m_F]^\C \cap \m_F^\C$ and 
it is equal to $1$ otherwise.
Note that the case $\ell_i = 2$
may occur only if $\g_F = {\goth f}_4$ or $\goth{sp}_n$ - see 
Table 1.
\endroster
\endproclaim
\demo{Proof} The first claim is an immediate consequence of Theorem 3.2 d) and 
the property of invariant complex structures on flag manifolds.\par
In order to prove  (3.12), let us consider 
a pair $\{\alpha_i, \alpha_i^d\}$ of CR dual roots in $R^+_F$ with $\alpha_i \neq \alpha^d_i$; 
by the previous remarks, 
there exist two integers $\ell_i$, $\ell^d_i$,which are
either $+1$ or $+2$, and two integers 
 $\epsilon_{\alpha_i}, \epsilon_{\alpha_i^d} = \pm 1$, so that 
$$E_{\alpha_i} + \epsilon_{\alpha_i} e^{2 \ell_i t}E_{- \alpha^d_i}\ , 
\ E_{\alpha^d_i} + \epsilon_{\alpha^d_i} e^{2 \ell^d_i t} E_{-\alpha_i} \in \m^{10}_F(t)$$
for any $t\neq 0$.\par
By direct inspection of Table 1, one can  check that the 
integers $\ell^d_i$, $\ell_i$ are  always equal. We claim that also $\epsilon_i = 
\epsilon_i^d$ for any CR dual pair $\{\alpha_i, \alpha^d_i\} \subset R^{(+)}_F$.\par
In fact,   by conjugation, it follows  that the following two vectors are in $\m^{01}_F(t)$ for any 
$t\neq 0$:
$$E_{\alpha_i} + \frac{1}{\epsilon_{\alpha^d_i} e^{2 \ell_i t}} E_{-\alpha^d_i}\ ,
\ E_{\alpha^d_i} + 
\frac{1}{\epsilon_{\alpha_i} e^{2 \ell_i t}}  E_{-\alpha_i} \in \m^{01}_F(t)\ .\tag 3.14$$
At this point, we recall that $\eta_0$ is a singular point for the $G$-action and that, 
by the structure theorems in \cite{HS} (see also \cite{AS}), 
the isotropy subalgebra $\g_{\eta_0}$ contains  the isotropy subalgebra $(\g_F)_{\eta_0}$
of the non-complex
singular $G_F$-orbit in $M$, which is a c.r.o.s.s.. In particular, 
one can check that  
$\dim_\R (\g_F)_{\eta_0} = \dim_\R \l_F + \dim_\C \m^{01}_F$. \par
On the other hand, by Lemma 3.5 (2), 
we have that $(\g_F)_{\eta_0} = \l_F + \g\cap \m^{01}_F(0)$ and hence that 
$$\dim_\R (\g\cap \m^{01}_F(0)) = \dim_\C \m^{01}_F(0)\ .\tag 3.15$$
Here, by $\m^{01}_F(0)$ we denote the subspace which is obtained 
from Table 1, by setting the value of the parameter $\lambda$ equal to 
 $\lambda(0) = e^{0} = 1$. Note that this subspace is {\it not} a 
Morimoto-Nagano subspace.\par
From (3.14), one can check that (3.15) occurs 
if and only if 
$$\epsilon_{\alpha^d_i} = \epsilon_{\alpha_i}\tag 3.16$$
for any pair of CR dual roots $\alpha_i, \alpha^d_i$. This proves the claim.\par
\medskip
In all the following, we will use the notation  $\epsilon_i = \epsilon_{\alpha_i} = 
\epsilon_{\alpha^d_i}$.\par
\medskip
By some straightforward computation, it follows that, for any $t\neq 0$, 
the elements $F_{\alpha_i}$, $F_{\alpha^d_i}$, 
$G_{\alpha_i}$ and $G_{\alpha^d_i}$ are equal to the following linear combinations
of  holomorphic and 
anti-holomorphic elements:
$$F_{\alpha_i} = 
\frac{1}{\sqrt{2}(1 -  e^{4 \ell_i t})}
\left\{\left[(E_{\alpha_i} + \epsilon_i e^{2\ell_i t}E_{-\alpha^d_i}) + 
\epsilon_i e^{2\ell_i t} (E_{\alpha^d_i} + \epsilon_i e^{2\ell_i t} E_{-\alpha_i})\right] + \right.$$
$$ + \left.
\left[-  e^{4 \ell_i t}
(E_{\alpha_i} + \frac{1}{\epsilon_i e^{2\ell_i t}} E_{-\alpha^d_i}) - 
\epsilon_i e^{2\ell_i t}(E_{\alpha^d_i} + 
\frac{1}{\epsilon_i e^{2\ell_i t}}  E_{-\alpha_i} )\right]\right\}\ ,$$
$$F_{\alpha^d_i} = 
\frac{1}{\sqrt{2}(1 -  e^{4 \ell_i t})}
\left\{\left[\epsilon_i e^{ 2\ell_i t} (E_{\alpha_i} + \epsilon_i e^{2\ell_i t}E_{-\alpha^d_i}) +
(E_{\alpha^d_i} + \epsilon_i e^{2\ell_i t} E_{-\alpha_i})\right] -\right.$$
$$ - \left.
\left[ e^{ 2\ell_i t} \epsilon_i
(E_{\alpha_i} + \frac{1}{\epsilon_i e^{2\ell_i t}} E_{-\alpha^d_i})
+  e^{4 \ell_i t} (E_{\alpha^d_i} + 
\frac{1}{\epsilon_i e^{2\ell_i t}}  E_{-\alpha_i} )\right]\right\}\ ,$$
$$G_{\alpha_i} = 
\frac{i}{\sqrt{2}(1 -  e^{4 \ell_i t})}
\left\{\left[(E_{\alpha_i} + \epsilon_i e^{2\ell_i t}E_{-\alpha^d_i})  - 
\epsilon_i e^{2\ell_i t} 
(E_{\alpha^d_i} + \epsilon_i e^{2\ell_i t} E_{-\alpha_i})\right] + \right.$$
$$ + \left.
\left[-  e^{4 \ell_i t}
(E_{\alpha_i} + \frac{1}{\epsilon_i e^{2\ell_i t}} E_{-\alpha^d_i}) + 
\epsilon_i e^{2\ell_i t} (E_{\alpha^d_i} + 
\frac{1}{\epsilon_i e^{2\ell_i t}}  E_{-\alpha_i} )\right]\right\}\ ,$$
$$G_{\alpha^d_i} = 
\frac{i}{\sqrt{2}(1 -  e^{4 \ell_i t})}
\left\{\left[- \epsilon_i e^{2\ell_i t} (E_{\alpha_i} + \epsilon_i e^{2\ell_i t}E_{-\alpha^d_i})  + 
(E_{\alpha^d_i} + \epsilon_i e^{2\ell_i t} E_{-\alpha_i})\right] + \right.$$
$$ + \left.
\left[ \epsilon_i e^{2\ell_i t}
(E_{\alpha_i} + \frac{1}{\epsilon_i e^{2\ell_i t}} E_{-\alpha^d_i}) -
 e^{4 \ell_i t}(E_{\alpha^d_i} + 
\frac{1}{\epsilon_i e^{2\ell_i t}}  E_{-\alpha_i} )\right]\right\}\ .$$
We then obtain that
$$J_t F_{\alpha_i} =  \frac{1 +  e^{4 \ell_i t}}
{1 -  e^{4 \ell_i t}} G_{\alpha_i} +
\frac{2 \epsilon_i e^{2\ell_i t}}
{1 -  e^{4 \ell_i t}} G_{\alpha^d_i}\ , $$
$$J_t F_{\alpha^d_i}  = 
\frac{ 2 \epsilon_i e^{ 2\ell_i t}}
{1 -  e^{4 \ell_i t}}  G_{\alpha_i} + 
\frac{1 +   e^{4 \ell_i t}}
{1 -  e^{4 \ell_i t}} G_{\alpha^d_i}\ .\tag 3.17$$
So, using the fact that $\epsilon_i^2 = 1$, 
we get  
$J_t F^+_i = \frac{1 +  e^{2\ell_i t}}{1 - e^{2\ell_i t}} G^+_i = - \coth(\ell_i t) G^+_i$ and 
$J_t F^-_{i} = \frac{1 -  e^{2\ell_i t}}{1 +  e^{2\ell_i t}} G^-_{i} = 
- \tanh(\ell_i t) G^-_i$.
The proof of (3.13) is similar.  It suffices to observe that for any $t\neq 0$
$$F^+_{i} = 
\frac{1}{\sqrt{2}(1 - e^{4\ell_i t})}
\left\{(1 + e^{2\ell_i t}) (E_{\alpha_i} +  e^{2\ell_i t}E_{- \alpha_i}) -\right.$$
$$\left.
- e^{2\ell_i t}(1 + e^{2\ell_i t})
(E_{\alpha_i} +  e^{-2\ell_i t} E_{-\alpha_i}) \right\}\ ,$$
$$G^+_{i} = 
\frac{i}{\sqrt{2}(1 - e^{4\ell_i  t})}
\left\{(1 - e^{2\ell_i t}) (E_{\alpha_i} +  e^{2\ell_i t}E_{- \alpha_i}) + 
 \right.$$
$$\left. +
e^{2\ell_i t}(1 - e^{2\ell_i t})
(E_{\alpha_i} +  e^{-2\ell_i t} E_{-\alpha_i}) \right\}\ ,$$
and hence that $J_t F^+_i = \frac{1 + e^{2\ell_i  t}}{1-e^{2\ell_i  t}} G^+_i = - \coth(\ell_i 
t) G^+_i$.\qed
\enddemo
\bigskip
\bigskip
\subhead 4. The algebraic representatives of the K\"ahler   and  Ricci form of a K-manifold
\endsubhead
\bigskip
In this section we give  a rigorous definition 
of the {\it algebraic representatives\/} of the K\"ahler form $\omega$ and the Ricci form $\rho$ 
of a K-manifold. We will also prove  Proposition 1.1.\par
\medskip
Indeed, we will give
the concept of 'algebraic representative'  for any bounded,
 closed 
2-form $\varpi$, which is  defined on  $M_{\text{reg}}$
and which is $G$-invariant and $J$-invariant.  
Clearly, 
$\omega|_{M_{\text{reg}}}$ and $\rho|_{M_{\text{reg}}}$ belong to this 
class of 2-forms.\par
\smallskip
Let  $\eta: \R \to M$ be an optimal transversal curve. Since $\g$ is semisimple, 
for any $G$-invariant 2-form $\varpi$  on $M_{\text{reg}}$
there exists a unique 
$\operatorname{ad}_\l$-invariant
element  $F_{\varpi, t}\in \Hom(\g,\g)$ such that:
$$\B(F_{\varpi, t}(X), Y) = \varpi_{\eta_t}(\hat X, \hat Y)\ ,\qquad X, Y \in \g\ ,\ t \neq 0\ .\tag4.1$$
If $\varpi$ is  also closed, 
we have that for any $X, Y, W\in \g$
$$0 = 3 d\varpi(\hat X, \hat Y, \hat W) = \varpi(\hat X, [\hat Y, 
\hat W]) +  \varpi(\hat Y, [\hat W, 
\hat X]) +  \varpi(\hat W, [\hat X, 
\hat Y])\ .$$
This implies that
$$F_{\varpi,t}([X,Y]), W)  = [F_{\varpi, t}(X), Y] + [X, F_{\varpi, t}(Y)]$$
i.e. $F_{\varpi, t}$ is a derivation of $\g$. Therefore, $F_{\varpi, t}$
is of the form
$$F_{\varpi, t} = \operatorname{ad}(Z_\varpi(t))\tag 4.2$$
for some $Z_\varpi(t) \in \g$ and 
$\varpi_{\eta_t}(\hat X, \hat Y) = \B([Z_\varpi(t), X], Y) = 
\B(Z_\varpi(t), [X,Y])$ . Note that since $F_{\varpi, t}$ is 
$\operatorname{ad}_\l$-invariant, then $Z_\varpi(t) \in C_\g(\l) = \z(\l) + \a$, where
$\a = C_\g(\l) \cap \l^\perp$.\par
We call the curve 
$$Z_\varpi: \R \to C_\g(\l) = \z(\l) + \a\ ,\tag 4.3$$
the {\it algebraic representative of the 2-form $\varpi$ along the 
optimal transversal curve $\eta$\/}.\par
 By definition, if the algebraic representative $Z_\varpi(t)$
is
given, it is possible to reconstruct 
 the values of $\varpi$ on any pair of vectors, which are
 tangent to
the regular orbits $G\cdot \eta_t$. Actually, since
for any  point $\eta_t \in M_{\text{reg}}$ we have that $J(T_{\eta_t} G) = T_{\eta_t} M$, 
it follows that 
 one can evaluate  $\varpi$  on {\it any\/} pair of 
vectors in $T_{\eta_t} M$  if the value
$\varpi_{\eta_t}(\hat Z_\D, J \hat Z_\D)$
is also given. However, in case 
$\varpi$ is a closed form,  the following Proposition shows
that  this last  value can 
be recovered from the first derivative of the function $Z_\varpi(t)$.\par
\bigskip
\proclaim{Proposition 4.1} Let $(M, J, g)$ be a K-manifold 
 acted on by the compact semisimple Lie group $G$. Let also $\eta_t = \exp(t i Z_\D)\cdot p_o$ 
be an optimal transversal curve
and  $Z_\varpi : \R \to \z(\l) + \a$
the algebraic representative of a bounded, $G$-invariant, $J$-invariant
closed 2-form $\varpi$ along $\eta$.  Then:\par
\roster
\item if $M$ is  a standard K-manifold or a 
non-standard KO-manifold (i.e. if either
 $\a = \R Z_\D$ or $\a = \goth{su_2}$ and $M$ is standard), 
then there exists an element $I_\varpi \in \z(\l)$ and 
a smooth function ${f}_\varpi: \R \to \R$ so that
$$Z_\varpi(t) = f_\varpi(t) Z_\D +  I_\varpi\ ;\tag 4.4$$
\item if $M$ is non-standard KE-manifold, then there exists 
a Cartan subalgebra $\t^\C \subset \l^\C + \a^\C$ and a root $\alpha$
of the corresponding root system, such  that  $Z_\D \in \R(iH_\alpha)$ 
and $\a = \R Z_\D + \R F_\alpha + \R G_\alpha$; furthermore
there exists an element $I_\varpi \in \z(\l)$, a real number $C_\varpi$ and 
a smooth function $f_\varpi: \R \to \R$ so that
$$Z_\varpi(t) = f_\varpi(t) Z_\D + \frac{C_\varpi}{\cosh(t)}
G_\alpha +  I_\varpi\ .\tag 4.4'$$
\endroster
Conversely,  if $Z_\varpi: \R \to C_\g(\l)$ is a curve in $C_\g(\l)$ of the form (4.4)
 or  (4.4'), 
then there exists a 
unique closed $J$-invariant, $G$-invariant 2-form 
$\varpi$  on $M_{\text{reg}}$, having $Z_\varpi(t)$ as  algebraic representative;  
such 2-form is the 
unique $J$- and  $G$-invariant form which 
verifies
$$\varpi_{\eta_t}(\hat V, \hat W) = \B(Z_\varpi(t), [V, W])\ ,
\qquad \varpi_{\eta_t}(J \hat Z_\D, \hat Z_\D)
= - f'_\varpi(t) \B(Z_\D,Z_\D)\ .\tag4.5$$
for any $V, W \in \m$ and any $\eta_t \in M_{\text{reg}}$.
\endproclaim
\demo{Proof} Let $\varpi$ be a closed 2-form which is $G$-invariant and $J$-invariant and let
$Z_\varpi(t)$ be the associated algebraic representative along $\eta$. Recall that $Z_\varpi(t)
\in \z(\l) + \a$. So,  if the action is ordinary (i.e. $\a = \R Z_\D$), $Z_\varpi(t)$
is of  the form
$$Z_\varpi(t) = f_\varpi(t) Z_\D +  I_\varpi(t)\ ,\tag 4.6$$
where the vector $I_\varpi(t) \in \z(\l)$  may depend on $t$. \par
In case the action of $G$ is extraordinary (that is 
 $\a =\goth{su}_2$)
by Lemma 2.2 in \cite{PS}, 
there exists  a Cartan subalgebra $\t^\C \subset \l^\C + \a^\C$, such that 
$\a^\C  = \C H_\alpha + \C E_\alpha
+ \C E_{-\alpha}$ for some root $\alpha$ of the corresponding root system. 
By the arguments in the proof of Theorem 3.7, this Cartan 
subalgebra can be always chosen in such a way that
$Z_\D \in \R(iH_\alpha)$ and hence that $\a = \R Z_\D + \R F_\alpha + 
\R G_\alpha$.\par
Then the function  $Z_\varpi(t)$ can be written as
$$Z_\varpi(t) = f_\varpi(t) Z_\D + g_\varpi(t) F_\alpha + h_\varpi(t) G_\alpha + I_\varpi(t)\tag4.6'$$
for some smooth real valued 
functions $f_\varpi$, $g_\varpi$ and $h_\varpi$ and  
some element $I_\varpi(t) \in \z(\l)$.\par
We now want  to show that, in case $M$ is a non-standard KE-manifold, 
then $g_\varpi(t) \equiv 0$ and 
that $h_\varpi(t) = \frac{C_\varpi}{\cosh(t)}$ for some constant $C_\varpi$.\par
In fact, observe that if   $Z_\varpi(t)$ is of the form (4.6') and if
$Z_\D$ is as listed in Table 1 for $\g_F = \goth{su}_2$, then 
$$\varpi_{\eta_t}(\hat Z_\D,\hat G_{\alpha}) = g_\varpi(t) \B(F_\alpha,[Z_\D, G_\alpha]) = 
- g_\varpi(t)\ ,$$
$$\varpi_{\eta_t}(\hat Z_\D, \hat F_\alpha) = h_\varpi(t) \B(G_\alpha, [Z_\D, F_\alpha]) = 
h_\varpi(t) \ .$$
Consider now the facts that $\varpi$ is closed, 
$\hat G_\alpha$ and $\hat Z_\D$ are holomorphic vector
fields and 
 $J \hat Z_\D|_{\eta_t} =  \eta'_t$. It follows that  
$g_\varpi$ verifies the following ordinary 
differential equation
$$\left.{d g_\varpi\over dt}\right|_{\eta_t}
= - \left.{d\over dt}\varpi(\hat Z_\D,\hat G_\alpha)\right|_{\eta_t} = 
-  \left.J\hat Z_\D \left(\varpi(\hat Z_\D,\hat G_\alpha)\right)\right|_{\eta_t} = $$
$$ = 
\left.\hat G_\alpha(\varpi(J\hat Z_\D,\hat Z_\D))\right|_{\eta_t} +
\left.\hat Z_\D(\varpi(\hat G_\alpha,J\hat Z_\D))\right|_{\eta_t} -
\varpi_{\eta_t}([J\hat Z_\D,\hat Z_\D],\hat G_\alpha) -$$
$$ -\left.
\varpi_{\eta_t}([\hat G_\alpha,J\hat Z_\D],\hat Z_\D) -
\varpi_{\eta_t}([\hat Z_\D,\hat G_\alpha],J\hat Z_\D)\right] = $$
$$ =  \varpi_{\eta_t}([\hat Z_\D,\hat G_\alpha],J\hat Z_\D) = 
- \varpi_{\eta_t}(\widehat{[Z_\D,G_\alpha]},J\hat Z_\D) = $$
$$ 
= - \varpi_{\eta_t}(\hat Z_\D, J\hat F_\alpha)  
= \coth(t)\varpi_{\eta_t}(\hat Z_\D, \hat G_\alpha) = 
- \coth(t)  g_\varpi(t)\ . 
\tag 4.7$$
We claim that this implies 
$$g_\varpi(t) \equiv 0\ .\tag 4.8$$
 In fact, if we assume that $g_\varpi(t)$ does not 
vanish identically, integrating the above equation, 
we have that
$g_\varpi(t) = \frac{C}{|\sinh(t)|}$ for some $C \neq 0$ and hence with 
a singularity at $t = 0$. But this 
contradicts the fact that $\varpi$ is a bounded 2-form. \par
With a similar argument, we have that $h_\varpi(t)$ verifies the differential equation 
$$\left.{d h_\varpi\over dt}\right|_{\eta_t} = 
-  \tanh(t) h_\varpi(t)\  ;$$
by 
integration this gives 
$$h_\varpi(t) = \frac{C_\varpi}{\cosh(t)}\tag 4.9$$
for some constant $C_\varpi$. \par
We show now that, in case $M$ is a  standard KE-manifold, then $Z_\varpi(t)$ 
is of the form (4.4). In fact, even if a priori $Z_\varpi(t)$ is of the form (4.6'), 
from Lemma 3.8 and the same arguments for proving
(4.7), we obtain that
$$\left.{d g_\varpi\over dt}\right|_{\eta_t}
= - \varpi_{\eta_t}(\hat Z_\D, J\hat F_\alpha)  
= - \varpi_{\eta_t}(\hat Z_\D, \hat G_\alpha) = 
g_\varpi(t)\ . 
\tag 4.10$$
This implies that $g_\varpi(t) = A e^{t}$ for some constant $A$. On the other hand,  
 if $A\neq 0$, it would follow that 
$\lim_{t\to \infty} |\varpi_{\eta_t}(\hat Z_\D, \hat G_\alpha)| = 
\lim_{t\to \infty} |g_{\varpi}(t)| = + \infty$, 
which is impossible since $\varpi_{\eta_t}(\hat Z_\D, \hat G_\alpha)$ is  bounded. Hence
$g_\varpi(t) \equiv 0$. \par
A similar argument proves that $h_\varpi(t) \equiv 0$.\par
\medskip
In order to conclude the proof, it remains
to show that in all cases the element $I_\varpi(t)$ is independent on $t$ and 
that $\varpi_{\eta_t}(J \hat Z_\D, Z_\D) = - f'_\varpi(t)
\B(Z_\D, Z_\D)$ for any $t$. 
We will prove these two facts 
 only for the case $\a \simeq \goth{sl}_2(\R)$ and $M$ non-standard, since the  proof 
in all other cases  is  similar. \par
Consider  two elements $V, W \in \g$. Since $\varpi$ is closed we have that 
$$0 = 3 d\varpi_{\eta_t}(J \hat Z_\D, \hat V, \hat W) 
 = $$
$$ = J \hat Z_\D(\varpi_{\eta_t}(\hat V, \hat W)) - \hat V(\varpi_{\eta_t}(J \hat Z_\D, \hat W))
+ W(\varpi_{\eta_t}(J \hat Z_\D, \hat V)) - $$
$$ - \varpi_{\eta_t}([J \hat Z_\D, \hat V], \hat W) + 
\varpi_{\eta_t}([J \hat Z_\D, \hat W], \hat V)  - \varpi_{\eta_t}([\hat V, \hat W], J \hat Z_\D) = $$
$$ = J \hat Z_\D|_{\eta_t}(\varpi(\hat V, \hat W)) - \varpi_{\eta_t}(J \hat Z_\D,[\hat V, \hat W]) = $$
$$ = \left.\frac{d}{dt}\left(\B(Z_\varpi, [V,W])\right)\right|_{t}
+ \varpi_{\eta_t}(J \hat Z_\D, \widehat{[V,W]})\ .\tag 4.11$$
On the other hand, we have the following orthogonal decomposition of the element $[V, W]$:
$$[V, W] = \frac{\B(Z_\D, [V, W])}{\B(Z_\D, Z_\D)} Z_\D -   
\B(F_\alpha, [V, W])  F_\alpha -  \B(G_\alpha, [V, W]) G_\alpha + $$
$$ +
[V, W]_{(\l + \a)^\perp} + [V, W]_\l\ ,$$
where $[V, W]_\l$ and $[V, W]_{(\l + \a)^\perp}$ are the 
orthogonal projections  of $[V, W]$
into $\l$  
and  $(\l + \a)^\perp$, respectively. Then 
$$\varpi_{\eta_t}(J \hat Z_\D, \widehat{[V,W]}) = 
\frac{\B(Z_\D, [V, W])}{\B(Z_\D, Z_\D)} \varpi_{\eta_t}(J \hat Z_\D, \hat Z_\D) - 
\B(F_\alpha, [V, W]) \varpi_{\eta_t}(J \hat Z_\D, \hat F_\alpha) - $$
$$ -
\B(G_\alpha, [V, W]) \varpi_{\eta_t}(J \hat Z_\D,\hat  G_\alpha) + 
\varpi_{\eta_t}(J \hat Z_\D, \widehat{[V, W]}_{(\l + \a)^\perp}) = $$
$$ = \frac{\B(Z_\D, [V, W])}{\B(Z_\D, Z_\D)} \varpi_{\eta_t}(J \hat Z_\D, \hat Z_\D) +
\B(F_\alpha, [V, W]) \varpi_{\eta_t}(\hat Z_\D, J\hat F_\alpha) + $$
$$ +
\B(G_\alpha, [V, W]) \varpi_{\eta_t}(\hat Z_\D, J \hat  G_\alpha) -
\varpi_{\eta_t}(\hat Z_\D, J \widehat{[V, W]}_{(\l + \a)^\perp}) = $$
$$ = \frac{\B(Z_\D, [V, W])}{\B(Z_\D, Z_\D)}  \varpi_{\eta_t}(J \hat Z_\D, \hat Z_\D) 
 +\B(G_\alpha, [V, W]) \frac{C_\varpi \tanh(t) }{\cosh(t)} -$$
$$
-
\B(Z_\varpi(t), \left[Z_\D, J_{\eta_t}([V, W]_{(\l + \a)^\perp})\right]) = $$
$$ = \B(\left\{ \frac{\varpi_{\eta_t}(J \hat Z_\D, \hat Z_\D)}{\B(Z_\D, Z_\D)} 
Z_\D - h'_\varpi(t)
 G_\alpha\right\}, [V, W])\ .$$
Therefore (4.11) becomes
$$\B\left(\left\{f'_\varpi(t)  + 
\frac{\varpi_{\eta_t}(J \hat Z_\D, \hat Z_\D)}{\B(Z_\D, Z_\D)} \right\} Z_\D 
+  \frac{d I_\varpi}{dt}, [V,W]\right) = 0\ .$$
Since $V, W$ are arbitrary 
 and $\frac{d I_\varpi}{dt} \in \z(l) \subset (Z_\D)^\perp$, 
 it implies  
$$f'_\varpi(t) =  - \frac{\varpi_{\eta_t}( J \hat Z_\D, \hat Z_\D)}{\B(Z_\D, Z_\D)}
\ ,\qquad
 \frac{d I_\varpi}{dt} \equiv 0\ ,$$
as we needed to prove. \qed
\enddemo
\bigskip
We conclude this section, with the following corollary which gives a geometric 
interpretation  of the optimal bases (see also \S 1).\par
\medskip
\proclaim{Corollary 4.2} Let $(M,J,g)$ be a K-manifold
and let $(F_i, G_i)$ be an optimal basis along an optimal transversal 
curve $\eta_t = \exp(t i Z) \cdot p_o$.  For any $\eta_t\in M_{\text{reg}}$, denote by
$\Cal F_t = (e_0, e_1, \dots, e_{n})_t$  the following 
holomorphic frame in  $T^\C_{\eta_t} M$: 
$$e_0 = \hat F_0|_{\eta_t} - i J \hat F_0|_{\eta_t} = \hat Z|_{\eta_t} - i J \hat Z|_{\eta_t}\ ,
\quad e_i = \hat F_i|_{\eta_t}  - i J \hat F|_{\eta_t} \ \ i\geq 1\ \ .$$
Then, 
\roster
\item if $M$ is a KO-manifold or a standard KE-manifold, then the holomorphic 
 frames $\Cal F_t$ are orthogonal w.r.t. any  $G$-invariant 
K\"ahler metric $g$ on $M$; 
\item if $M$ is a non-standard KE-manifold, then 
the holomorphic 
 frames $\Cal F_t$ are orthogonal w.r.t. any  $G$-invariant 
K\"ahler metric $g$ on $M$, whose associated algebraic representative 
$Z_\omega(t)$  has  vanishing coefficient   $C_\omega = 0$ (see Proposition 4.1 for the definition
 of $C_\omega$)
\endroster
\endproclaim
\demo{Proof} It is a direct consequence of  definitions and 
Proposition 4.1.\qed
\enddemo
\bigskip
\bigskip
\subhead 5. The Ricci tensor  of  a K-manifold
\endsubhead
\bigskip
From the results of \S 4,  the Ricci form $\rho$ can be
completely recovered from 
the algebraic representative $Z_\rho(t)$ along an optimal transversal
curve $\eta_t$.  
On the other 
hand, using a few  known properties of flag manifolds, the reader can check that 
the curve  $Z_\rho(t) \in \z(\l) + \a$ is uniquely  determined
by the 1-parameter family of quadratic forms $Q^r$ on  $\m$ given by 
$$Q^r_{t}: \m \to \R\ ,\qquad 
Q^r_t(E) = r_{\eta_t}(\hat E, \hat E)\qquad \left( = - \rho_{\eta_t}(\hat E, \hat E) = 
- \B(Z_\rho(t), [E, J_t E]) \right)  
\ .$$
Since $\m$ corresponds to the subspace $\D_{\eta_t} \subset T_{\eta_t} G\cdot \eta_t$, 
this means that {\it for any 
K\"ahler metric $\omega$, the corresponding the Ricci tensor $r$ is uniquely determined by 
its
restrictions $r|_{\D_t \times \D_t}$
on the  holomorphic tangent spaces $\D_t$ of the regular orbits $G\cdot \eta_t$\/}.\par
\medskip
The expression for the restrictions $r|_{\D_t \times \D_t}$ 
in terms of the algebraic representative $Z_\omega(t)$ of the 
 K\"ahler form $\omega$ is given in the following Theorem.   
\par
\bigskip
\proclaim{Theorem 5.1} Let $(M, J, g)$ be a K-manifold
and  $\eta_t = \exp(ti Z_\D)\cdot p_o$ be an optimal 
transversal curve. Using the same notation of \S 3, 
let also 
$(F_i, G_i) = 
(F_0, F^{\pm}_k, G^{\pm}_k, F'_j, G'_j)$ be an optimal basis for $\R Z_\D + \m$; finally, 
 for any $1\leq j\leq N_F$ let $\ell_j$ be the integer
which appear in (3.12) for  the expression of
$J_t F_i$ and  for any $N_F+1 \leq k \leq n-1$ let  $\beta_k$ be the root so that 
$F_k = F_{\beta_k}$.\par
Then, for any $\eta_t \in M_{\text{reg}}$ and for any element 
$E \in \m$
$$\rho_{\eta_t}(\hat E , J \hat E) =  A_{E}(t) \left\{ \frac{1}{2} h'(t)  -
 \sum_{i = 1}^{N_F} \tanh^{(-1)^{i+1}}(\ell_i t) \ell_i +
\sum_{j = N_F+1}^{n-1} \beta_j(iZ_\D)\right\} + B_E(t)
\tag 5.1$$
where 
$$h(t) = \log(\omega^n(\hat F_0, J \hat F_0, \hat F_1, J \hat F_1, \dots , J F_{n-1})|_{\eta_t}) 
\ ,\tag 5.2$$
$$A_{E}(t) =\frac{\B([E, J_t E], Z_\D)}{\B(Z_\D, Z_\D)}\ , \tag 5.3$$
$$B_E(t) =  
- \sum_{i = 1}^{N_F} \tanh^{(-1)^{i+1}}(\ell_i t)\B([E,J_t E]_{\l + \m} , [F_i, G_i]_{\l + \m} )
+ $$
$$ + \sum_{j = N_F+1}^{n-1} \B\left(i H_{\beta_j}, [E,J_t E]_{\z(\l)}\right)\ ,
\tag 5.4$$
and where,
for any $X\in \g$, we denote by 
$X_{\l+ \m}$ (resp. $X_{\z(\l)}$) the projection parallel to $(\l + \m)^\perp = \R Z_\D$ 
(resp. to $\z(\l)^\perp$) of $X$ into 
$\l + \m$ (resp. into $\z(\l)$). 
\endproclaim
\medskip 
\demo{Proof} 
Let  
 $J_t$ be  the complex structure on $\m$ induced by the complex structure 
$J$ of $M$. For any $E\in \m$
and any point $\eta_t$, we  may clearly write that 
$\rho_{\eta_t}(\hat E, J \hat E) =  \rho_{\eta_t}(\hat E , \widehat{J_t E} )$ 
and hence,  by   Koszul's formula 
 (see \cite{Ko}, \cite{Be}), 
$$\rho_{\eta_t}(\hat E , J \hat E ) =
\frac{1}{2} \frac {\left(\Cal L_{J \widehat{[ E, J_t E]}}
\omega^n\right)_{\eta_t}(\hat F_0, J \hat F_0, \hat F_1, J \hat F_1, \dots , J F_{n-1})}
{\omega^n_{\eta_t}(\hat F_0, J \hat F_0, \hat F_1, J \hat F_1, \dots , J F_{n-1})}
\tag 5.5$$
(note that the definition we adopt here for the Ricci form $\rho$ is opposite 
in sign to the definition used in \cite{Be}).\par
Recall that for any $Y \in \g$, 
we 
may write 
$$\hat Y|_{\eta(t)} = \sum_{i \geq  0} \lambda_i \hat F_i
|_{\eta(t)} + \sum_{i \geq 1}\mu_i J \hat F_i|_{\eta(t)}\ ,$$
where
$$\lambda_i = \frac{\B(Y, F_i)}{\B(F_i,F_i)} \ ,\qquad
\mu_i = \frac{\B(Y, J_t F_i)}{\B(J_t F_i,J_tF_i)}\ .$$
Hence, for any $i$
$$[J \widehat{[ E, J_t E]}, \hat F_i]_{\eta_t} = - J  \widehat{[[ E, J_t E], F_i]}_{\eta_t} = $$
$$ = 
- \sum_{j \geq  0}  \frac{\B([ [ E, J_t E], F_i], F_j)}{\B(F_j,F_j)} J \hat F_j
|_{\eta(t)} + \sum_{j \geq 1}\frac{\B([[ E, J_t E], F_i]
, J_t F_j)}{\B(J_t F_j,J_tF_j)}\hat F_j|_{\eta(t)} = $$
$$ = - \sum_{j \geq  0}  \frac{\B([ E, J_t E], [F_i, F_j])}{\B(F_j,F_j)} J \hat F_j
|_{\eta(t)} + \sum_{j \geq 1}\frac{\B([ E, J_t E], [F_i
, J_t F_j])}{\B(J_t F_j,J_tF_j)}\hat F_j|_{\eta(t)}\ , \tag 5.6$$
$$[J \widehat{[ E, J_t E]}, J \hat F_i]_{\eta_t} = 
\widehat{[ [ E, J_t E], F_i]}_{\eta_t} = $$
$$ = 
 \sum_{j \geq  0}  \frac{\B( [ E, J_t E],[ F_i, F_j])}{\B(F_j,F_j)} \hat F_j
|_{\eta(t)} + \sum_{j \geq 1}\frac{\B([ E, J_t E], [F_i
, J_t F_j])}{\B(J_t F_j,J_tF_j)} J\hat F_j|_{\eta(t)}\ .\tag 5.7$$
Therefore, if we denote 
$h(t) = \log(\omega^n(\hat F_0, J \hat F_0, \hat F_1, J \hat F_1, \dots, J F_{n-1} )|_{\eta_t})$, 
then, after some straightforward computations, (5.5) becomes
$$\rho_{\eta_t}(\hat E , \widehat{J_t E} ) =
\frac{1}{2}J \widehat{[ E, J_t E]}(h)|_{\eta_t} - 
\sum_{i \geq 1}^{n-1}
\frac{\B([ E, J_t E], [F_i
, J_t F_i])}{\B(J_t F_i,J_tF_i)}
\ .
\tag 5.8$$
We claim  that 
$$J \widehat{[ E, J_t E]}(h)|_{\eta_t} = A_E(t) h'_t\ .\tag 5.9$$
In fact, 
 for any 
$X\in \g$
$$\hat X(\omega(\hat F_0, J\hat F_0, \dots, J \hat F_{n-1})|_{\eta_t} = $$
$$ =
- \omega_{\eta_t}(\widehat{[X, F_0]}, J \hat F_0, \dots, J \hat F_{n-1})  -
\omega( F_0, J \widehat{[X, F_0]}, \dots, J \hat F_{n-1}) - \dots = 0\ . \tag 5.10$$
On the other hand, 
$$ J \widehat{[ E, J_t E]}|_{\eta_t} = 
\frac{\B([ E, J_t E], Z_\D)}{\B(Z_\D, Z_\D)} J \hat Z_\D|_{\eta_t} 
+ J \hat X_{\eta_t} =
A_E(t) J \hat Z_\D|_{\eta_t} 
+ \widehat{J_t X}_{\eta_t} \tag 5.11$$
 for some  some $X \in \m$. 
From (5.11) and (5.10) and the fact that $J \hat Z_\D|_{\eta_t} = \eta'_t$, 
we immediately obtain (5.9).\par
\medskip
Let us now prove that
$$\sum_{i \geq 1}^{n-1}
\frac{\B([ E, J_t E], [F_i
, J_t F_i])}{\B(J_t F_i,J_tF_i)}=  A_E\left\{ 
\sum_{i = 1}^{N_F} \tanh^{(-1)^{i+1}}(\ell_i t) \ell_i - 
\sum_{j = N_F+1}^{n-1} \beta_i(iZ_\D)\right\} - B_E\tag 5.12$$
First of all, observe that from definitions, for any 
$1 \leq k \leq N_F$ we have that, for any case of Table 1,  when $\alpha_k \neq \alpha_k^d$, 
$$\B(Z_\D, [F_k, G_k]) = \frac{1}{2} \B(Z_\D, [F_{\alpha_k} + (-1)^{k+1} 
\epsilon_k F_{\alpha_k^d}, G_{\alpha_k} + (-1)^{k+1} 
\epsilon_k G_{\alpha_k^d}]) = $$
$$ = \frac{i}{2} \B(Z_\D, H_{\alpha_k} + H_{\alpha^d_k}) = \ell_k\ ,\tag 5.13$$
and, when $\alpha_k = \alpha_k^d$, 
$$\B(Z_\D, [F_k, G_k]) = \B(Z_\D, [F_{\alpha_k}, G_{\alpha_k}]) = 
\B(Z_\D, i H_{\alpha_k}) =  \ell_k\ .\tag 5.13'$$
Similarly, for any $N_F + 1 \leq j \leq n-1$
$$\B(Z_\D, [F_j, G_j]) = \B(Z_\D, i H_{\beta_j}) = \beta_j(i Z_\D)\ .\tag 5.14$$
So, using (5.13), (5.13'), (5.14) and the fact that $\B(F_i, F_i) = \B(G_i, G_i) = -1$ for any 
$1 \leq i \leq n-1$, we obtain that 
 for  $1 \leq k \leq N_F$, 
$$\frac{\B([ E, J_t E], [F_k
, J_t F_k])}{\B(J_t F_k,J_tF_k)} = 
\tanh^{(-1)^{k+1}}(\ell_k t) \left(
\B(Z_\D, [F_k, G_k])
\frac{\B([ E, J_t E],Z_\D)}{\B(Z_\D, Z_\D)}  + \right.$$
$$ + \left. \phantom{\frac{\B([ E, J_t E],Z_\D)}{\B(Z_\D, Z_\D)}}
\B([ E, J_t E]_{\l+\m}, [F_k, G_k]_{\l + \m})\right) = $$
$$ = \tanh^{(-1)^{k+1}}(\ell_k t)  \left[A_E(t)\ell_k +
\B([ E, J_t E]_{\l+\m}, [F_k, G_k]_{\l + \m})\right]\ ,\tag 5.15$$
and
for any $N_F +1 \leq j \leq N$
$$\frac{\B([ E, J_t E], [F_j
, J_t F_j])}{\B(J_t F_j,J_tF_j)} = $$
$$ = 
- \frac{\B([ E, J_t E],Z_\D)}{\B(Z_\D, Z_\D)} \B(Z_\D, 
[F_j
, G_j]) - \B([ E, J_t E]_{\l+ \m}[F_j
, G_j]_{\l+ \m}) =$$
$$ = - A_E \beta_j(iZ_\D) - \B(i H_{\beta_j}, [ E, J_t E]_{\z(\l)})\ . \tag 5.16$$
 From (5.15) and 
(5.16), we immediately obtain (5.12) and from (5.8) this concludes the proof. \qed
\enddemo
\bigskip
The expressions for the functions $A_E(t)$ and $B_E(t)$
simplify considerably if one assumes that $E$ is an element of the optimal basis. Such expressions
are given  
in the following conclusive proposition.\par
\medskip
\proclaim{Proposition 5.2} Let $(F_i, G_i)$ be an optimal 
basis along an optimal transversal curve $\eta_t$ of a K-manifold
$M$. For any $1\leq i \leq N_F$, let $\ell_i$ be as in Theorem 5.1 and 
 denote  by
 $\{\alpha_i, \alpha^d_i\}\subset R^{(+)}_F$ 
the pair of CR-dual roots, such that  
$F_i = \frac{1}{\sqrt{2}}(F_{\alpha_i} \pm \epsilon_i F_{\alpha^d_i})$ or 
$F_i = F_{\alpha_i}$, in case $\alpha_i = \alpha^d_i$; 
also, for any $N_F +1 \leq j \leq n-1$,   denote
 by $\beta_j \in R'_+$  the root such that 
$F_j = F_{\beta_j}$. Finally,  let  
$A_E(t)$ and $B_E(t)$ be as
defined in Theorem 5.1 and let us denote by 
$$Z^\kappa = \sum_{k = N_F + 1}^{n-1} i H_{\beta_k}\ . \tag 5.17$$
\roster
\item If $E = F_i$ 
for some  $1 \leq i \leq N_F$, then 
$$A_{F_i}(t) =  - \frac{\ell_i \tanh^{(-1)^{i}}(\ell_i t) }{
\B(Z_\D, Z_\D)}\ ;\tag5.18$$
and 
$$B_{F_i}(t)  =  - \frac{\ell_i \tanh^{(-1)^i}(\ell_i t)}{\B(Z_\D, Z_\D)} \B(Z^\kappa, Z_\D) + $$
$$ + 
\tanh^{(-1)^{i}}(\ell_i t)\left(
\sum_{j = 1}^{N_F} \tanh^{(-1)^{j+1}}(\ell_j t)\B([F_i,G_i]_{\l + \m} , [F_j, G_j]_{\l + \m} )\right)\ ,
\tag 5.19$$
\item If $E = F_i$ for some  $N_F +1 \leq i \leq n-1$, then
$$A_{F_i}(t) = \frac{\B(Z_\D, i H_{\beta_i})}{\B(Z_\D, Z_\D)} \ ,  
\quad
B_{F_i}(t) = \B(Z^\kappa, i H_{\beta_i})\ .\tag 5.20$$ 
\endroster
\endproclaim
\demo{Proof} Formulae (5.18) and $(5.19)$ are immediate consequences of
definitions and of (5.13), (5.13') and (5.14). Formula (5.20)
can be checked using the fact that 
$[F_{\beta_i}, J_t F_{\beta_i}] = [F_{\beta_i}, G_{\beta_i}] = i H_{\beta_i}$ for any 
$N_F + 1 \leq i \leq n-1$, from
 properties of the Lie brackets $[F_i, G_i]$, with $1 \leq i \leq N_F$, which 
can be derived from Table 1, and from the fact that $\R Z_\D \subset [\m', \m']^\perp$. 
\qed
\enddemo

\enddocument
\bye